\documentclass{compositio} 
\usepackage{amsmath,amscd}
\begin{document}

%%%%%%%%%%%%%%%%%%%%%%%%%%%%%%%%%%%%%%%%%%%%%%%%%%%%%%%%%%%%%%%%%%%%%%
%% START OF MACRO DEFINITIONS
%%%%%%%%%%%%%%%%%%%%%%%%%%%%%%%%%%%%%%%%%%%%%%%%%%%%%%%%%%%%%%%%%%%%%%

%%%%%%%%%%%%%%%%%%%%%%%%%%%%%%%%%%%%%%%%%%%%%%%%%%%%%%%%%%%%%%%%%%%%%%
% Theorem environments

\newtheorem{theorem}{Theorem}
\newtheorem{lemma}[theorem]{Lemma}
\newtheorem{conjecture}[theorem]{Conjecture}
\newtheorem{question}[theorem]{Question}
\newtheorem{proposition}[theorem]{Proposition}
\newtheorem{corollary}[theorem]{Corollary}

\theoremstyle{definition}
% The * surpresses numbering
\newtheorem*{definition}{Definition}
\newtheorem{example}[theorem]{Example}

\theoremstyle{remark}
\newtheorem{remark}[theorem]{Remark}

\newcommand\claim[2]{\par\vspace{1ex minus 0.5ex}\noindent%
\textbf{Claim #1}.\enspace\emph{#2}.\par\noindent\ignorespaces}

%%%%%%%%%%%%%%%%%%%%%%%%%%%%%%%%%%%%%%%%%%%%%%%%%%%%%%%%%%%%%%%%%%%%%%

%%%%%%%% Set Up Environment for Notation %%%%%%%%%%%%%%
% This is currently set to allow quite wide items to be defined
\newenvironment{notation}[0]{%
  \begin{list}%
    {}%
    {\setlength{\itemindent}{0pt}
     \setlength{\labelwidth}{4\parindent}
     \setlength{\labelsep}{\parindent}
     \setlength{\leftmargin}{5\parindent}
     \setlength{\itemsep}{0pt}
     }%
   }%
  {\end{list}}

%%%%%%%% Set Up Environment for Parts in Theorems %%%%%%%%%%%%%%
\newenvironment{parts}[0]{%
  \begin{list}{}%
    {\setlength{\itemindent}{0pt}
     \setlength{\labelwidth}{1.5\parindent}
     \setlength{\labelsep}{.5\parindent}
     \setlength{\leftmargin}{2\parindent}
     \setlength{\itemsep}{0pt}
     }%
   }%
  {\end{list}}
% Use \Part{(a)}, instead of \item[(a)], to ensure upright font
\newcommand{\Part}[1]{\item[\upshape#1]}

%%%%%%%%%%%%%%%%%%
% Greek Alphabet %
%%%%%%%%%%%%%%%%%%
\renewcommand{\a}{\alpha}
\renewcommand{\b}{\beta}
\newcommand{\g}{\gamma}
\renewcommand{\d}{\delta}
\newcommand{\e}{\epsilon}
\newcommand{\f}{\phi}
\renewcommand{\l}{\lambda}
\renewcommand{\k}{\kappa}
\newcommand{\lhat}{\hat\lambda}
\newcommand{\m}{\mu}
\newcommand{\bfmu}{{\boldsymbol{\mu}}}
\renewcommand{\o}{\omega}
\renewcommand{\r}{\rho}
\newcommand{\rbar}{{\bar\rho}}
\newcommand{\s}{\sigma}
\newcommand{\sbar}{{\bar\sigma}}
\renewcommand{\t}{\tau}
\newcommand{\z}{\zeta}

\newcommand{\D}{\Delta}
\newcommand{\G}{\Gamma}
\newcommand{\F}{\Phi}

%%%%%%%%%%%%%%%%%%%%
% Fraktur Alphabet %
%%%%%%%%%%%%%%%%%%%%
\newcommand{\ga}{{\mathfrak{a}}}
\newcommand{\gb}{{\mathfrak{b}}}
\newcommand{\gn}{{\mathfrak{n}}}
\newcommand{\gp}{{\mathfrak{p}}}
\newcommand{\gP}{{\mathfrak{P}}}
\newcommand{\gq}{{\mathfrak{q}}}

%%%%%%%%%%%%%%%%%%%
% Barred Alphabet %
%%%%%%%%%%%%%%%%%%%
\newcommand{\Abar}{{\bar A}}
\newcommand{\Ebar}{{\bar E}}
\newcommand{\Pbar}{{\bar P}}
\newcommand{\Sbar}{{\bar S}}
\newcommand{\Tbar}{{\bar T}}
\newcommand{\ybar}{{\bar y}}
\newcommand{\phibar}{{\bar\phi}}

%%%%%%%%%%%%%%%%%%%%%%%%%
% Calligraphic Alphabet %
%%%%%%%%%%%%%%%%%%%%%%%%%
\newcommand{\Acal}{{\mathcal A}}
\newcommand{\Bcal}{{\mathcal B}}
\newcommand{\Ccal}{{\mathcal C}}
\newcommand{\Dcal}{{\mathcal D}}
\newcommand{\Ecal}{{\mathcal E}}
\newcommand{\Fcal}{{\mathcal F}}
\newcommand{\Gcal}{{\mathcal G}}
\newcommand{\Hcal}{{\mathcal H}}
\newcommand{\Ical}{{\mathcal I}}
\newcommand{\Jcal}{{\mathcal J}}
\newcommand{\Kcal}{{\mathcal K}}
\newcommand{\Lcal}{{\mathcal L}}
\newcommand{\Mcal}{{\mathcal M}}
\newcommand{\Ncal}{{\mathcal N}}
\newcommand{\Ocal}{{\mathcal O}}
\newcommand{\Pcal}{{\mathcal P}}
\newcommand{\Qcal}{{\mathcal Q}}
\newcommand{\Rcal}{{\mathcal R}}
\newcommand{\Scal}{{\mathcal S}}
\newcommand{\Tcal}{{\mathcal T}}
\newcommand{\Ucal}{{\mathcal U}}
\newcommand{\Vcal}{{\mathcal V}}
\newcommand{\Wcal}{{\mathcal W}}
\newcommand{\Xcal}{{\mathcal X}}
\newcommand{\Ycal}{{\mathcal Y}}
\newcommand{\Zcal}{{\mathcal Z}}

%%%%%%%%%%%%%%%%%%%%%%%%%%%%
% Blackboard Bold Alphabet %
%%%%%%%%%%%%%%%%%%%%%%%%%%%%
\renewcommand{\AA}{\mathbb{A}}
\newcommand{\BB}{\mathbb{B}}
\newcommand{\CC}{\mathbb{C}}
\newcommand{\FF}{\mathbb{F}}
\newcommand{\GG}{\mathbb{G}}
\newcommand{\PP}{\mathbb{P}}
\newcommand{\QQ}{\mathbb{Q}}
\newcommand{\RR}{\mathbb{R}}
\newcommand{\ZZ}{\mathbb{Z}}

%%%%%%%%%%%%%%%%%%%%%%%%%%
% Boldface Math Alphabet %
%%%%%%%%%%%%%%%%%%%%%%%%%%
\newcommand{\bfa}{{\mathbf a}}
\newcommand{\bfb}{{\mathbf b}}
\newcommand{\bfc}{{\mathbf c}}
\newcommand{\bfe}{{\mathbf e}}
\newcommand{\bff}{{\mathbf f}}
\newcommand{\bfg}{{\mathbf g}}
\newcommand{\bfp}{{\mathbf p}}
\newcommand{\bfr}{{\mathbf r}}
\newcommand{\bfs}{{\mathbf s}}
\newcommand{\bft}{{\mathbf t}}
\newcommand{\bfu}{{\mathbf u}}
\newcommand{\bfv}{{\mathbf v}}
\newcommand{\bfw}{{\mathbf w}}
\newcommand{\bfx}{{\mathbf x}}
\newcommand{\bfy}{{\mathbf y}}
\newcommand{\bfz}{{\mathbf z}}
\newcommand{\bfA}{{\mathbf A}}
\newcommand{\bfF}{{\mathbf F}}
\newcommand{\bfB}{{\mathbf B}}
\newcommand{\bfG}{{\mathbf G}}
\newcommand{\bfI}{{\mathbf I}}
\newcommand{\bfM}{{\mathbf M}}
\newcommand{\bfzero}{{\boldsymbol{0}}}

%%%%%%%%%%%%%%%%%%%%%%%%%%%%%%
% Miscellaneous New Commands %
%%%%%%%%%%%%%%%%%%%%%%%%%%%%%%
\newcommand{\Aut}{\operatorname{Aut}}
\newcommand{\CM}{\operatorname{CM}}   % CM(O) = set of t with O_t = O
\newcommand{\Disc}{\operatorname{Disc}}
\newcommand{\Div}{\operatorname{Div}}
\newcommand{\Ell}{\operatorname{Ell}}   % Ell(O) = set of E with CM by O
\newcommand{\End}{\operatorname{End}}
\newcommand{\Fbar}{{\bar{F}}}
\newcommand{\Gal}{\operatorname{Gal}}
\newcommand{\GL}{\operatorname{GL}}
\newcommand{\Index}{\operatorname{Index}}
\newcommand{\Image}{\operatorname{Image}}
\newcommand{\liftable}{{\textup{liftable}}}
\newcommand{\hhat}{{\hat h}}
\newcommand{\Ker}{{\operatorname{ker}}}
\newcommand{\Lift}{\operatorname{Lift}}
\newcommand{\MOD}[1]{~(\textup{mod}~#1)}
\newcommand{\Norm}{{\operatorname{\mathsf{N}}}}
\newcommand{\notdivide}{\nmid}
\newcommand{\normalsubgroup}{\triangleleft}
\newcommand{\odd}{{\operatorname{odd}}}
\newcommand{\onto}{\twoheadrightarrow}
\newcommand{\ord}{\operatorname{ord}}
\newcommand{\Pic}{\operatorname{Pic}}
\newcommand{\Prob}{\operatorname{Prob}}
\newcommand{\Qbar}{{\bar{\QQ}}}
\newcommand{\rank}{\operatorname{rank}}
\newcommand{\Resultant}{\operatorname{Resultant}}
\renewcommand{\setminus}{\smallsetminus}
\newcommand{\Span}{\operatorname{Span}}
\newcommand{\tors}{{\textup{tors}}}
\newcommand{\Trace}{\operatorname{Trace}}
\newcommand{\UHP}{{\mathfrak{h}}}    % Upper half plane
\newcommand{\<}{\langle}
\renewcommand{\>}{\rangle}

\newcommand{\longhookrightarrow}{\lhook\joinrel\longrightarrow}
\newcommand{\longonto}{\relbar\joinrel\twoheadrightarrow}

%%%%%%%%%%%%%%%%%%%%%%%%%%%%%%%%%%%%%%%%%%%%%%%%%%%%%%%%%%%%%%%%%%%%%%
%% END OF MACRO DEFINITIONS
%%%%%%%%%%%%%%%%%%%%%%%%%%%%%%%%%%%%%%%%%%%%%%%%%%%%%%%%%%%%%%%%%%%%%%

\title
[On the independence of Heegner points]
{On the independence of Heegner points associated to
     distinct quadratic imaginary fields}

%% First Author
\author{Michael Rosen}
\email{mrosen@math.brown.edu}
\address{Mathematics Department\\Box 1917\\
         Brown University\\Providence, RI 02912 USA}
%% Second Author
\author{Joseph H. Silverman}
\email{jhs@math.brown.edu}
\address{Mathematics Department\\Box 1917\\
         Brown University\\Providence, RI 02912 USA}

\thanks{The second author's research supported by NSA grant H98230-04-1-0064}

\begin{abstract}
Let $E/\QQ$ be an elliptic curve with a fixed modular parametrization
$\F_E:X_0(N)\to{E}$ and let~$P_1,\ldots,P_r\in{E}(\Qbar)$ be Heegner points
attached to the rings of integers of distinct quadratic imaginary
field~$k_1,\ldots,k_r$. We prove that if the odd parts of
the class numbers of~$k_1,\ldots,k_r$ are larger than a constant~$C=C(E,\F_E)$
depending only on~$E$ and~$\F_E$, then the points~$P_1,\ldots,P_r$ 
are independent in~$E(\Qbar)/E_\tors$. We also discuss a possible 
application to the elliptic curve discrete logarithm problem.
\end{abstract}

%%%%%%%%%%%%%%%%%%%%%%%%%%%%%%%%%%%%%%%%%%%%%%%%%%%%%%%%%%%%%%%%%%%%%%%%
%% NON-TEX ABSTRACT FOR ARXIV
%
%% Let E/Q be an elliptic curve with a fixed modular parametrization
%% F : X_0(N) --> E and let P_1,...,P_r be Heegner points on E
%% attached to the rings of integers of distinct quadratic imaginary
%% field k_1,...,k_r. We prove that if the odd parts of
%% the class numbers of k_1,...,k_r are larger than a constant C=C(E,F)
%% depending only on E and F, then the points P_1,...,P_r 
%% are independent in E/(torsion). We also discuss a possible 
%% application to the elliptic curve discrete logarithm problem.
%
%%%%%%%%%%%%%%%%%%%%%%%%%%%%%%%%%%%%%%%%%%%%%%%%%%%%%%%%%%%%%%%%%%%%%%%%

\maketitle

%%%%%%%%%%%%%%%%%%%%%%%%%%%%%%%%%%%%%%%%%%%%%%%%%%%%%%%%%%%%%%%%%%%%%%%%
% Section: Introduction
%%%%%%%%%%%%%%%%%%%%%%%%%%%%%%%%%%%%%%%%%%%%%%%%%%%%%%%%%%%%%%%%%%%%%%%%
\section*{Introduction}

The theory of Heegner points provides a fundamental method for
creating algebraic points on modular curves and on the elliptic curves
that they parametrize. The work of Wiles
et.al.~\cite{BreuilConradDiamondTaylor01,ConradDiamondTaylor99,TaylorWiles95,Wiles95}
says that every elliptic curve~$E/\QQ$ of conductor~$N$ admits a modular
parametrization $\F_E:X_0(N)\to{E}$, so in particular there is a
theory of Heegner points on elliptic curves defined over~$\QQ$.
Heegner points appear prominently in the work of Gross, Kohnen, and
Zagier~\cite{GrossKohnenZagier87,GrossZagier86} on the
Birch-Swinnerton-Dyer conjecture and in the work of
Kolyvagin~\cite{Kolyvagin88b,Kolyvagin88a} on Mordell-Weil ranks and
Shafarevich-Tate finiteness.  (A nice survey of this material may be
found in~\cite{DarmonCBMS04}.)
\par
Of particular importance in the work of Kolyvagin and others is the
construction of Euler systems of Heegner points.
(See~\cite{Kolyvagin90} for a general formulation.)  These are
collections of Heegner points~$(P_n)$ defined over a tower of ring
class fields lying over a \emph{single} quadratic imaginary field and
satisfying various trace and Galois compatibility conditions. In this
paper we consider the orthogonal problem of collections of Heegner
points~$(P_n)$ defined over class fields of \emph{different} quadratic
imaginary fields.
\par
Our main theorem says that under a fairly mild class number condition,
such a set of Heegner points on~$E$ corresponding to distinct
quadratic imaginary fields has maximal rank in~$E(\Qbar)/E_\tors$. We
briefly state the result here and refer the reader to
Sections~\ref{section:heegnerpointsonX}
and~\ref{section:heegnerpointsonE} for definitions and to
Section~\ref{section:indepheegnerpts} for a more precise statement.

\begin{theorem}
\label{theorem:heegindepintro}
Let~$E/\QQ$ be an elliptic curve
\textup{(}with a given modular parametrization\textup{)}.
There is a constant~$C=C(E)$ so that the following is true.
\par
Let~$k_1,\ldots,k_r$ be distinct quadratic imaginary fields whose
class numbers satisfy
\[
  h(k_i)^\odd \ge C
  \qquad\text{for all $1\le i\le r$,}
\]
where~$h^\odd$ denotes the odd part of the integer~$h$.
Let~$P_1,\ldots,P_r\in{E}(\Qbar)$ be Heegner points associated to
\textup{(}the ring of integers of\textup{)}~$k_1,\ldots,k_r$,
respectively.  Then
\[
  \text{$P_1,\dots,P_r$ are independent in $E(\Qbar)/E_\tors$.}
\]
\end{theorem}

\begin{remark}
Theorem~\ref{theorem:heegindepintro} is possibly not surprising, and
indeed the statement may be true with no (or a much weaker) class
number hypothesis. On the other hand, since the compositum of the
quadratic imaginary fields~$k_1,\ldots,k_r$ may have degree as small
as~$r$, as opposed to the maximal value of~$2^r$, there seems no
obvious reason why the associated Heegner points must be completely
independent.  Thus elementary considerations might lead to an estimate
of the form
\[
  2^{\rank(\ZZ P_1+\cdots+\ZZ P_r)} \ge [k_1\cdots k_r : \QQ],
\]
but the proof of the stronger statement given in
Theorem~\ref{theorem:heegindepintro} requires a blend of
class field theory, Galois theory, linear algebra
(over~$\ZZ/n\ZZ$), and Serre's theorem on the image of Galois
in~$\Aut(E_\tors)$.
\end{remark}

\begin{remark}
Other authors have considered the behavior of Heegner points
associated to different quadratic imaginary fields. In particular, we
mention the fundamental work of Gross, Kohnen and
Zagier~\cite{GrossKohnenZagier87}. In the notation of
Theorem~\ref{theorem:heegindepintro}, the Heegner
point~$P_i\in{E}(\Qbar)$ is defined over the Hilbert class field~$K_i$
of~$k_i$. We can obtain points defined over~$\QQ$ by taking the trace,
\[
  Q_i = \Trace_{K_i/\QQ}(P_i) \in E(\QQ).
\]
Gross, Kohnen, and Zagier~[\emph{op.cit.}] compute the canonical
height pairing~$\langle{Q_i},{Q_j}\rangle$ of these points and prove
that~$Q_1,\ldots,Q_r$ generate a subgroup of~$E(\QQ)$ of rank at
most~$1$. More precisely, they show that
\[
  \rank(\ZZ Q_1+\cdots+\ZZ Q_r) = \begin{cases}
    1&\text{if }\ord_{s=1}L(E/\QQ,s)=1,\\
    0&\text{otherwise.}
  \end{cases}
\]
This is in accordance with the predictions of the
Birch-Swinnerton-Dyer conjecture.
\end{remark}

Aside from its intrinsic interest, the independence result of
Theorem~\ref{theorem:heegindepintro} has a (negative) application to
the elliptic curve discrete logarithm problem (ECDLP). If the theorem
were false and Heegner points had a tendency to be dependent, then
potentially there would be an algorithm to solve the ECDLP on elliptic
curves with small coefficients by using Deuring lifts and Heegner
points. We briefly sketch the idea in Section~\ref{section:ecdlp}
and refer the reader to~\cite{RosenSilverman05b}
for a more detailed description.

Finally, in Section~\ref{section:2partclassgroup}, we make some brief
remaraks and raise a question concerning the distribution of quadratic
imaginary fields whose class numbers have bounded odd parts.

\begin{acknowledgements}
The authors would like to thank Gary Cornell and Siman Wong for
their helpful suggestions.
\end{acknowledgements}

%%%%%%%%%%%%%%%%%%%%%%%%%%%%%%%%%%%%%%%%%%%%%%%%%%%%%%%%%%%%%%%%%%%%%%%%
% Section: Heegner points on $X_0(N)$
%%%%%%%%%%%%%%%%%%%%%%%%%%%%%%%%%%%%%%%%%%%%%%%%%%%%%%%%%%%%%%%%%%%%%%%%
\section{Heegner points on $X_0(N)$}
\label{section:heegnerpointsonX}

In this section we briefly review the theory of Heegner points on the
modular curve~$X_0(N)$ and in the next section we discuss Heegner
points on elliptic curves.  We refer the reader
to~\cite[\S\S3.1,3,3]{DarmonCBMS04} and~\cite{Gross84} for further
details.  \par Recall that the noncuspidal points of the modular
curve~$X_0(N)$ classify isomorphism classes of triples~$(A,A',\f)$
consisting of two elliptic curves $A$ and~$A'$ and an
isogeny~\text{$\f:A\to{A}'$} whose kernel is cyclic of order~$N$.
Heegner points are associated to orders in quadratic imaginary fields,
so we set
\begin{notation}
\item[$k/\QQ$]
a quadratic imaginary field,
\item[$\Ocal_k$]
the ring of integers of $k$,
\item[$\Ocal$]
an order in~$\Ocal_k$.
\end{notation}
Every order has the form~\text{$\Ocal=\ZZ+c\Ocal_k$} for a unique
integer~$c\ge1$ called the \emph{conductor of~$\Ocal$}. The
discriminant of~$\Ocal$ is given by
\[
  \Disc(\Ocal) = c^2\Disc(\Ocal_k),
  \quad\text{where $c^2=(\Ocal_k:\Ocal)$.}
\]
\par
In order to describe the Heegner points on~$X_0$,
we follow that notation used in~\cite{DarmonCBMS04} and define:
\par\noindent
\begin{tabular}{cl}
$\Pic(\Ocal)$
    &\parbox[t]{.8\hsize}{%
    the \emph{Picard group} (or \emph{class group}) of $\Ocal$,
    defined to be the group of isomorphism classes of rank~1
    projective $\Ocal$-modules. If $\Ocal=\Ocal_k$, then $\Pic(\Ocal)$
    is the usual ideal class group of~$\Ocal_k$.
    } \\
$\Ell^{(N)}(\Ocal)$
    &\parbox[t]{.8\hsize}{%
    the set of isomorphism classes of triples~$(A,A',\f)$ such
    that~$A$ and~$A'$ are elliptic curves satisfying 
    \vspace*{1\jot}\hfill\break\hspace*{4em}
    $\End(A)\cong\End(A')\cong\Ocal$ 
    \vspace*{1\jot}\hfill\break
    and $\f:A\to{A'}$ is an isogeny with $\Ker(\f)\cong\ZZ/N\ZZ$.
    } \\
\end{tabular}
\par\noindent
\begin{tabular}{cl}
$\CM^{(N)}(\Ocal)$
    &\parbox[t]{.8\hsize}{%
    the set of points in the noncuspidal part of~$X_0(N)$
    corresponding to the triples~$(A,A',\f)$ in $\Ell^{(N)}(\Ocal)$.
    } \\
\end{tabular}

\vspace{2\jot minus 1\jot}

\begin{definition}
We will generally identify without further comment  the two sets
\[
  \Ell^{(N)}(\Ocal)\longleftrightarrow \CM^{(N)}(\Ocal).
\]
The points in either set are called \emph{Heegner points} of~$X_0(N)$.
\end{definition}

It is clearly important to determine conditions on~$\Ocal$ that
ensure that there exist nontrivial Heegner points.

\begin{proposition}
Assume that the discriminant of~$\Ocal$ is prime to~$N$.  Then the set
$\Ell^{(N)}(\Ocal)$ is nonempty if and only if every prime
dividing~$N$ is split in~$k$. 
\end{proposition}
\begin{proof}
See~\cite[Proposition~3.8]{DarmonCBMS04}, or see~\cite[\S3]{Gross84}
for a stronger statement in which it is only required that~$N$ be
relatively prime to the conductor of~$\Ocal$..
\end{proof}

\begin{definition}
We say that~$\Ocal$ satisfies the \emph{Heegner condition for~$N$} if
the following two conditions are satisfied:
\begin{itemize}
\Part{(1)}
$\gcd\bigl(\Disc(\Ocal),N\bigr)=1$.
\Part{(2)}
Every prime dividing~$N$ is split in~$k$.
\end{itemize}
We say that~$k$ satisfies the Heegner condition for~$N$ if its ring of
integers satisfies the condition.
\end{definition}

Let~$\Ocal$ be an order in~$k$.  Class field theory associates
to~$\Ocal$ a finite abelian extension~$K_\Ocal/k$, called the
\emph{ring class field of~$k$ attached to~$\Ocal$.} The
extension~$K_\Ocal/k$ is unramified outside the primes dividing the
conductor of~$\Ocal$, and the Artin reciprocity map gives an
isomorphism
\begin{equation}
  \label{equation:artinrecip}
  (\;\cdot\;,K_\Ocal/k) : \Pic(\Ocal)
  \xrightarrow{\;\sim\;} \Gal(K_\Ocal/k).
\end{equation}
In particular, if~$\bar\gp\in\Pic(\Ocal)$ corresponds to a prime ideal
of~$k$ that does not divide~$\Disc(\Ocal)$,
then~$(\bar\gp,K_\Ocal/k)$ is the inverse of the Frobenius element
at~$\gp$.

\begin{theorem}
\label{theorem:heegnerptonX0}
Let~$\Ocal$ be an order in~$k$ that satisfies the Heegner condition for~$N$.
\begin{parts}
\Part{(a)}
The points in~$\CM^{(N)}(\Ocal)$ are defined over~$K_\Ocal$, i.e., 
\[
  \CM^{(N)}(\Ocal) \subset X_0(N)(K_\Ocal).
\] 
\Part{(b)}
The points of~$\Ell^{(N)}(\Ocal)$ are in one-to-one correspondence with
the set of pairs
\[
  \bigl\{(\gn,\bar\ga) : 
  \text{$\bar\ga\in\Pic(\Ocal)$, $\gn$ is a proper $\Ocal$-ideal,
             $\Ocal/\gn\cong\ZZ/N\ZZ$} \bigr\}.
\]
The correspondence is given explicitly by associating to a
pair~$(\gn,\bar\ga)$ the cyclic $N$-isogeny
\[
  \CC/\ga \xrightarrow{\;z\to z\;} \CC/\ga\gn^{-1}.
\]
\Part{(c)}
There is a natural action of~$\Pic(\Ocal)$ on~$\Ell^{(N)}(\Ocal)$
\textup{(}and thus also on~$\CM^{(N)}(\Ocal)$\textup{)} which we
denote by~$\star$. In terms of pairs~$(\gn,\bar\ga)$, it is given by
the formula
\[
  \bar\gb\star(\gn,\bar\ga) = (\gn,\bar\ga\bar\gb).
\]
\Part{(d)}
The~$\star$-action is compatible with the action of Galois via the
reciprocity map in the sense that 
\[
  y^{(\bar\gb,K_\Ocal/k)} = \bar\gb^{-1}\star y
  \qquad\text{for all~$y\in\CM^{(N)}(\Ocal)$.}
%%  \bar\gb^{-1}\star(\gn,\bar\ga) = (\gn,\bar\ga\bar\gb^{-1}).
\]
\end{parts}
\end{theorem}
\begin{proof}
See~\cite[Chapter~3]{DarmonCBMS04} or~\cite{Gross84}.
\end{proof}

For our purposes, the importance of
Theorem~\ref{theorem:heegnerptonX0} is that it allows us to conclude
that every point in~$\CM^{(N)}(\Ocal)$ generates a large extension
of~$k$, as in the following result.

\begin{corollary}
\label{corollary:heegnerptdegX0}
Let~$\Ocal$ be an order in~$k$ that satisfies the Heegner condition
for~$N$ and let $y\in\CM^{(N)}(\Ocal)$. Then
\[
  k(y) = K_\Ocal.
\]
\end{corollary}
\begin{proof}
We know from Theorem~\ref{theorem:heegnerptonX0}(a) that~$y$ is
defined over~$K_\Ocal$. Further, 
if we identify~$y$ with a pair~$(\gn,\bar\ga)$ as in
Theorem~\ref{theorem:heegnerptonX0}(b), then~(c) and~(d)
tell us that the full set of Galois conjugates of~$y$ is given
by
\begin{align*}
  \bigl\{y^\s : \s\in\Gal(K_\Ocal/k) \bigr\}
  &= \bigl\{\bar\gb\star y : \bar\gb\in\Pic(\Ocal) \bigr\} \\
  &= \bigl\{ (\gn,\bar\ga\bar\gb^{-1}) : \bar\gb\in\Pic(\Ocal) \bigr\} \\
  &= \bigl\{ (\gn,\bar\gb) : \bar\gb\in\Pic(\Ocal) \bigr\} 
\end{align*}
The points~$(\gn,\bar\gb)$ are distinct for distinct~$\gb\in\Pic(\Ocal)$,
so we see that
\begin{equation}
  \label{equation:Kkkykyspic}
  [K_\Ocal:k]
  \ge \bigl[k(y):k\bigr]
  = \#\bigl\{y^\s : \s\in\Gal(K_\Ocal/k) \bigr\} \ge \#\Pic(\Ocal).
\end{equation}
Class field theory~\eqref{equation:artinrecip}
tells us that~\text{$\#\Pic(\Ocal)=[K_\Ocal:k]$}.
Hence all of the inequalities in~\eqref{equation:Kkkykyspic}
are equalities, which proves that~$k(y)=K_\Ocal$.
\end{proof}

%%%%%%%%%%%%%%%%%%%%%%%%%%%%%%%%%%%%%%%%%%%%%%%%%%%%%%%%%%%%%%%%%%%%%%%%
% Section: Heegner points on elliptic curves
%%%%%%%%%%%%%%%%%%%%%%%%%%%%%%%%%%%%%%%%%%%%%%%%%%%%%%%%%%%%%%%%%%%%%%%%
\section{Heegner points on elliptic curves}
\label{section:heegnerpointsonE}

Let~$E/\QQ$ be an elliptic curve of conductor~$N$.
The theorem of Wiles et.al.~\cite{BreuilConradDiamondTaylor01,%
ConradDiamondTaylor99,%
TaylorWiles95,%
Wiles95}
says that there exists a modular parametrization
\[
  \F_E : X_0(N) \longrightarrow E.
\]
The map~$\F_E:X_0(N)\to{E}$ is a finite covering defined over~$\QQ$.

\begin{definition}
Let~$k$ be a quadratic imaginary field and let~$\Ocal$ be an order in~$k$
that satisfies the Heegner condition for~$N$. The set of \emph{Heegner
points of~$E$} (associated to~$\Ocal$) is the set
\[
  \bigl\{\F_E(y) : y\in\CM^{(N)}(\Ocal) \bigr\}.
\]
\end{definition}

The action of~$\Pic(\Ocal)$ and~$\Gal(K_\Ocal/K)$
on~$\CM^{(N)}(\Ocal)$ as described in Theorem~\ref{theorem:heegnerptonX0}
translates directly into analogous actions on Heegner points on~$E$,
see~\cite[Theorems~3.6,~3.7]{DarmonCBMS04}.  All that we will require
is the following elementary consequence.

%% \begin{theorem}
%% \label{theorem:galoisheegnerpt}
%% With notation as above, let~$\Ocal$ be an order in~$k$
%% and let~$\tau\in\CM^{(N)}(\Ocal)$.
%% \begin{parts}
%% \Part{(a)}
%% The Heegner point~$P_\tau=\F_E(\tau)\in E(\CC)$ is defined over the
%% ring class field~$K_\Ocal$.
%% \Part{(b)}
%% \textup{(Shimura reciprocity)}
%% The action of~$\Gal(K_\Ocal/k)$ on~$P_\tau$ is given,
%% via the reciprocity map, by the following formula:
%% \[
%%   \F_E(\tau)^{(\bar\ga,K_\Ocal/k)} = \F_E(\bar\ga^{-1}\star\tau) 
%%   \qquad\text{for $\bar\ga\in\Pic(\Ocal)$.}
%% \]
%% \end{parts}
%% \end{theorem}
%% \begin{proof}
%% See~\cite[Theorems~3.6,~3.7]{DarmonCBMS04}.
%% \end{proof}

\begin{proposition}
\label{proposition:heegnerptdegE}
Let~$\Ocal$ be an order in~$k$ that satsifies the Heegner condition
for~$N$, let~$y\in\CM^{(N)}(\Ocal)$, and let~$P_y=\F_E(y)$ be the
associated Heegner point. Then
\[
  \bigl[k(P_y):k\bigr] \ge \frac{[K_\Ocal:k]}{\deg\F_E}.
\]
\end{proposition}
\begin{proof}
To ease notation, let
\[
  d=\deg\F_E,\qquad
  n=  [K_\Ocal:k],\qquad
  m=  [k(P_y):k\bigr].
\]
From Corollary~\ref{corollary:heegnerptdegX0} we know
that~$k(y)=[K_\Ocal:k]$, so~$y$ has exactly~$n$ Galois conjugates,
say~$y_1,\dots,y_n$.  Further, Galois acts transitively
on the collection of points~$y_1,\dots,y_n$, so it acts transitively on their
images~$\F_E(y_1),\F_E(y_2),\dots,\F_E(y_n)$. Since~$\F_E$ is at most
$d$-to-1, it follows that~$\F(y)$ has at least~$n/d$ distinct
conjugates, and hence that~\text{$m\ge{n/d}$}.
\end{proof}

The ring class field~$K_\Ocal$ is an abelian extension of~$k$.  It is
not abelian extension of~$\QQ$, but it is a Galois extension,
and there is an exact sequence
\[
  1 \longrightarrow \Gal(K_\Ocal/k)
    \longrightarrow \Gal(K_\Ocal/\QQ)
    \longrightarrow \Gal(k/\QQ)
    \longrightarrow  1.
\]
The elements in~$\Gal(K_\Ocal/\QQ)$ with nontrivial image
in~$\Gal(k/\QQ)$ are called~\emph{reflections}.  The
structure of~$\Gal(K_\Ocal/\QQ)$ and its action on Heegner points
is described in the following proposition.

\begin{proposition}
Let~$\Ocal$ be an order in the quadratic imaginary field~$k$,
let~$K_\Ocal$ be the associated ring class field, and
let~$\rho\in\Gal(K_\Ocal/\QQ)$ be a reflection.
\begin{parts}
\Part{(a)}
The extension~$K_\Ocal/\QQ$ is Galois with group equal to a
generalized dihedral group. More precisely, 
\[
  \rho \s = \s^{-1}\rho.
  \qquad\text{for all $\s\in\Gal(K_\Ocal/k)$.}
\]
\Part{(b)}
Let $w(E/\QQ)$ be the sign of the functional equation of~$E/\QQ$,
let~$y\in\CM^{(N)}(\Ocal)$, and let~$P_\tau=\F_E(\t)$ be the
associated Heegner point on~$E$.  Then there exists
a~$\s\in\Gal(K_\Ocal/k)$, depending on~$P_y$ and~$\rho$, so
that
\[
  P_y^\rho = -w(E/\QQ)P_y^\s
  \pmod{ E(K_\Ocal)_\tors}.
\]
\end{parts}
\end{proposition}
\begin{proof}
See \cite[Proposition~3.11]{DarmonCBMS04} or~\cite{Gross84}.
\end{proof}

%%%%%%%%%%%%%%%%%%%%%%%%%%%%%%%%%%%%%%%%%%%%%%%%%%%%%%%%%%%%%%%%%%%%%%%%
% Section: A linear algebra estimate
%%%%%%%%%%%%%%%%%%%%%%%%%%%%%%%%%%%%%%%%%%%%%%%%%%%%%%%%%%%%%%%%%%%%%%%%
\section{A linear algebra estimate}
\label{section:linearalgebra}

The intuition behind the following proposition is that a large
subgroup of~$\GL_2(\ZZ/n\ZZ)$ cannot act in an abelian fashion on a
large subgroup of~$(\ZZ/n\ZZ)^2$. The key for our application is to
quantify this statement in such a way that it is uniform with respect
to~$n$.

\begin{proposition}
\label{proposition:orderindexbd}
Suppose that the following quantities are given:
\begin{notation}
\item[$n$]
a positive integer.
\item[$V$]
a free~$\ZZ/n\ZZ$ module of rank~$2$.
\item[$\G$]
a subgroup of $\Aut(V)\cong\GL_2(\ZZ/n\ZZ)$.
\item[$W$]
a $\G$-invariant $\ZZ/n\ZZ$-submodule of~$V$, i.e., $\G{W}=W$.
\end{notation}
Let
\[
  I(\G) = \bigl(\Aut(V):\G\bigr) 
  = \text{the index of $\G$ in~$\Aut(V)$.}
\]
Suppose that the action of~$\G$ on~$W$ is \emph{abelian} in the
sense that~$\G|_W$ is an abelian subgroup of~$\Aut(W)$.
Then
\begin{equation}
  \label{equation:WllIG4}
  |W| \le I(\G)^3.
\end{equation}
\end{proposition}

\begin{remark}
For our purposes it suffices to know that~$|W|$ is bounded in terms
of~$I(\G)$, independent of~$n$, but it is an interesting question to
ask to what extent the inequality~\eqref{equation:WllIG4} might be
improved.  A more elaborate argument (which we omit) gives
exponent~$2$, and one might hope for an estimate of the form
\text{$|W|\ll{I}(\G)^{1+\e}$}.  However, the following example shows
that an exponent of at least~$\frac{4}{3}$ is necessary.
\end{remark}

\begin{example}
Let~$\ell$ be a prime, let~$n=\ell^{2k}$, and let~$W=V$. Then one
easily checks that the action of
\[
  \G = \bigl\{ A\in\GL_2(\ZZ/\ell^{2k}\ZZ) : 
         \text{$A\bmod\ell^k$ is diagonal} \bigr\} 
\]
on~$(\ZZ/\ell^{2k}\ZZ)^2$ is abelian. Then
\begin{align*}
  |\G| &= \bigl|(\ZZ/\ell^k\ZZ)^*\bigr|\cdot\bigl|M_2(\ZZ/\ell^k\ZZ)\bigr|
    = \ell^{5k}(1-\ell^{-1}), \\
  \bigl|\GL(\ZZ/\ell^{2k}\ZZ)\bigr| &= \ell^{8k}(1-\ell^{-1})(1-\ell^{-2}), \\
  I(\G) &= \ell^{3k}(1-\ell^{-2}).
\end{align*}
Since~$|W|=\ell^{4k}$, this yields
\[
  \frac{\log|W|}{\log I(\G)} \xrightarrow[k\to\infty]{} \frac{4}{3}.
\]
\end{example}

\begin{proof}[Proof of Proposition \text{\ref{proposition:orderindexbd}}]
By the standard structure theorem on modules over
PIDs~\cite[Theorem~VI.2.7]{LangAlgebra} we can find a basis for~$V$ so
that
\[
  V \cong \frac{\ZZ}{n\ZZ}\times\frac{\ZZ}{n\ZZ}
  \qquad\text{and}\qquad
  W \cong \frac{m_1\ZZ}{n\ZZ}\times\frac{m_2\ZZ}{n\ZZ}
  \quad\text{with $m_1|m_2|n$.}
\]
We begin by doing the case that $n=\ell^e$ is a power of a prime, so
\[
  V \cong \frac{\ZZ}{\ell^e\ZZ}\times\frac{\ZZ}{\ell^e\ZZ}
  \quad\text{and}\quad
  W \cong \frac{\ell^i\ZZ}{\ell^e\ZZ}\times\frac{\ell^{i+j}\ZZ}{\ell^e\ZZ}
  \quad\text{with $i,j\ge0$ and $i+j\le{e}$.}
\]
This allows us to identify~$\G$ with a subgroup
of~$\GL_2(\ZZ/\ell^e\ZZ)$.  Further, the condition that~$\G{W}=W$
implies that every
matrix~$\left(\begin{smallmatrix}a&b\\c&d\\\end{smallmatrix}\right)\in\G$
satisfies~$c\equiv0\pmod{\ell^j}$.
In terms of the classical modular groups,
the condition~$\G{W}=W$ is equivalent to the requirement
\begin{equation}
  \label{equation:GinG0elljGL2}
  \G \subset \Image
    \bigl(\G_0(\ell^j)\longrightarrow\GL_2(\ZZ/\ell^e\ZZ) \bigr).
\end{equation}
\par
It may happen that~\eqref{equation:GinG0elljGL2} is true for a larger
value of~$j$, so we define
\[
  J = \min\Bigl\{e,\max\left\{\ord_\ell(c) : 
     \left(\begin{smallmatrix} a&b\\c&d\\\end{smallmatrix}\right) \in \G
   \right\} \Bigr\}
\]
In other words,~$J$ is the largest integer less than or equal to~$e$
such that every matrix in~$\G$ is congruent
to~$\left(\begin{smallmatrix} *&*\\0&*\\\end{smallmatrix}\right)$
modulo~$\ell^J$. 
\par
In particular, we can find a matrix
\[
  A = \begin{pmatrix} a&b\\ p^Jc&d\\ \end{pmatrix} \in \G
  \qquad\text{with $c\not\equiv0\pmod{\ell}$.}
\]
(Notice that if~$J=e$, we can simply take~$c=1$.)  We fix the
matrix~$A$ and consider another matrix~$B\in\G$. The definition of~$J$
tells us that~$B$ has the form
\[
  B = \begin{pmatrix} \a&\b\\ p^J\g&\d\\ \end{pmatrix}.
\]
If we were working over a field, we might hope that the condition
\text{$AB|_W=BA|_W$} implies that~$B$ has the form~$xI+yA$ for some
scalars~$x$ and~$y$.  This is not quite true in our case, but we will
now prove the validity of a similar statement up to a carefully chosen
power of~$\ell$.
\par
We consider first the case that
\[
  e > i+J
\]
and apply the assumption that~$A$ and~$B$ commute in their action
on~$W$.  Note that this is not the same as saying that~\text{$AB=BA$}
in~$\GL_2(\ZZ/\ell^e\ZZ)$. We obtain the correct statement by
requiring that~\text{$BA-AB$} kills a basis of~$W$. Thus
\[
  (BA-AB)\begin{pmatrix}\ell^i&0\\0&\ell^{i+j}\\\end{pmatrix}
  \equiv \begin{pmatrix}0&0\\0&0\\\end{pmatrix} \pmod{\ell^e}.
\]
Multiplying this out and doing some algebra yields
\begin{multline}
  \label{equation:BA-AB=0}
  \begin{pmatrix}
  \ell^{i+J}(c\b-b\g) & \ell^{i+j}\bigl(b(\a-\d)-(a-d)\b\bigr) \\
  \ell^{i+J}\bigl((a-d)\g-c(\a-\d)\bigr) & -\ell^{i+j+J}(c\b-b\g) \\
  \end{pmatrix}  \\
  \equiv \begin{pmatrix}0&0\\0&0\\\end{pmatrix} \pmod{\ell^e}.
\end{multline}
We use the congruences in the matrix equation~\eqref{equation:BA-AB=0}
to compute
\begin{align*}
  \ell^{i+J} \bigl((c\a-a\g)I+\g A\bigr)
  &= \ell^{i+J} \begin{pmatrix}
      c\a & b\g \\ \ell^J c \g & c\a - \g(a-d) \\
      \end{pmatrix} \\
  &\equiv \ell^{i+J} \begin{pmatrix}
      c\a & c\b \\ \ell^J c \g & c\d \\
      \end{pmatrix}  \pmod{\ell^e} \\
  &= \ell^{i+J} c B.
\end{align*}
Using the fact that~$\ell\notdivide{c}$, we have shown that 
every~$B\in\G$  satisfies
\[
  B \equiv xI+yA \pmod{\ell^{e-i-J}}
  \qquad\text{for some $x,y\in\ZZ/\ell^{e-i-J}\ZZ$,}
\]
i.e., every~$B\in\G$ has the form
\[
  B = xI+yA+\ell^{e-i-J}Z
  \quad\text{with}\quad\begin{cases}
     x,y\in\ZZ/\ell^{e-i-J}\ZZ\quad\text{and}\\
     Z\in M_2(\ZZ/\ell^{i+J}\ZZ). \\
   \end{cases}
\]
Of course, the group~$\G$ will not contain all of these matrices
(e.g., we cannot have $\ell|x$ and~$\ell|y$), but
in any case we obtain an upper bound
\[
  |\G| \le (\text{\# of $(x,y)$})\cdot(\text{\# of $Z$}) \\
       = \ell^{2(e-i-J)} \cdot \ell^{4(i+J)}\\
       = \ell^{2e+2i+2J}.
\]
The order of~$\GL_2(\ZZ/\ell^e\ZZ)$ is well known, so we obtain
a lower bound for the index
\begin{multline}
  \label{equation:indexbd1}
  I(\G) = \frac{\bigl|\GL_2(\ZZ/\ell^e\ZZ)\bigr|}{|\G|}
  \ge \frac{\ell^{4e}(1-\ell^{-1})(1-\ell^{-2})}{\ell^{2e+2i+2J}} \\
  = \ell^{2e-2i-2J}(1-\ell^{-1})^2(1+\ell^{-1}).
\end{multline}
This estimate is helpful provided that~$J$ is not too large. However, 
in the case that~$J$ is large, we can instead use the fact that~$\G$
is contained in the image of~$\G_0(\ell^J)$ to estimate
\[
  |\G| \le \left|\Big\{ 
    \left(\begin{smallmatrix}a&b\\c&d\\\end{smallmatrix}\right)
    \in\GL_2(\ZZ/\ell^e\ZZ) : \ell^J | c \Big\}\right|
  = \ell^{4e-J}(1-\ell^{-1})^2,
\]
so
\begin{equation}
  \label{equation:indexbd2}
  I(\G) 
  \ge \frac{\ell^{4e}(1-\ell^{-1})(1-\ell^{-2})}{\ell^{4e-J}(1-\ell^{-1})^2}
  = \ell^J \frac{1+\ell^{-1}}{1-\ell^{-1}}.
\end{equation}
Multiplying~\eqref{equation:indexbd1} by the square
of~\eqref{equation:indexbd2} yields
\[
  I(\G)^3 \ge \ell^{2e-2i}(1+\ell^{-1})^3,
\]
which proves that~$I(\G)^3$ is larger than
\[
  |W| = \ell^{2e-2i-j}.
\]
\par
Next we consider the case that 
\begin{equation}
  \label{equation:iJgee}
  i+J\ge e.
\end{equation}
The commutativity relation~\eqref{equation:BA-AB=0} then gives little
information, but the fact that~$\G$ is contained in the image
of~$\G_0(\ell^J)$ still gives the lower bound~\eqref{equation:indexbd2},
which we use in the weaker form~\text{$I(\G)\ge\ell^J$}. Combining this
with the assumption~\eqref{equation:iJgee} yields
\[
  I(\G)^2 \ge \ell^{2J} \ge \ell^{2e-2i} \ge \ell^{2e-2i-j} = |W|,
\]
which is stronger than the desired result.  This completes the proof
of the proposition in the case that~$n$ is a power of a prime.
\par
Finally, suppose that~$n$ is arbitrary. Let
\[
  V_\ell=V\otimes\ZZ_\ell=\text{$\ell$-primary part of $V$},
\]
and similarly let $W_\ell=W\otimes\ZZ_\ell$.  Then
\[
  V=\bigoplus_{\ell|n}{V}_\ell
  \qquad\text{and}\qquad
  W=\bigoplus_{\ell|n}{W}_\ell
\]
by the Chinese remainder theorem. Further, we have
\[
  \Aut(V)=\bigoplus_{\ell|n} \Aut(V_\ell) 
  \qquad\text{and}\qquad
  \G=\bigoplus_{\ell|n}\G_\ell\quad 
    \text{with $\G_\ell=\Image\Bigl(\G\to\Aut(V_\ell)\Bigr)$.}
\]
Applying the $\ell$-primary case to this direct sum decomposition yields
\[
  |W| = \prod_{\ell|n} |W_\ell|
  \le \prod_{\ell|n} I(\G_\ell)^3 = I(\G)^3,
\]
which completes the proof of Proposition~\ref{proposition:orderindexbd}.
\end{proof}

%%%%%%%%%%%%%%%%%%%%%%%%%%%%%%%%%%%%%%%%%%%%%%%%%%%%%%%%%%%%%%%%%%%%%%%%
% Section: Multiples of points and abelian extensions
%%%%%%%%%%%%%%%%%%%%%%%%%%%%%%%%%%%%%%%%%%%%%%%%%%%%%%%%%%%%%%%%%%%%%%%%
\section{Multiples of points and abelian extensions}
\label{section:independenceresult}

Our eventual goal is to prove the independence of points that are
defined over fields that are ``sufficiently large, sufficiently
disjoint, and sufficiently abelian.'' The content of the next theorem
is to quantify the meaning of the word ``sufficiently'' in this
statement.  Its proof combines the elementary linear algebra result
from Section~\ref{section:linearalgebra} with Serre's deep theorem on
the image of Galois.

\begin{theorem}
\label{theorem:abstractindepresult}
Let $F/\QQ$ be a number field, let~$E/F$ be an elliptic curve that
does not have complex multiplication, and let~$d\ge1$. There is an
integer~$M=M(E/F,d)$ so that for any field~$k/F$ and
any~$P\in{E}(\Fbar)$ satisfying
\[
  [k:F]\le d
  \text{\qquad and\qquad $k(P)/k$ is abelian,}
\]
the following estimate is true:
\[
  \bigl[k(P):k\bigr] \enspace\text{divides}\enspace M\bigl[k(nP):k\bigr] 
  \qquad\text{for all $n\ge1$.}
\]
\end{theorem}
\begin{proof}
Fix an integer $n\ge1$. Writing
\[
  \bigl[k(P):k\bigr] = \bigl[k(P):k(nP)\bigr] \bigl[k(nP):k\bigr],
\]
it suffices to find a bound of the form
\[
  \bigl[k(P):k(nP)\bigr] \le C = C(E/F,d),
\]
since then we can take~$M$ to equal the least common multiple of the
integers less then~$C$. 
\par
Consider the following set of points in~$E$,
\begin{equation}
  S(P,n) = \label{equation:SPnPtsPt}
    \bigl\{ P^{\t\s}-P^\t : 
       \s\in\Gal\bigl(k(P)/k(nP)\bigr),\ 
%%       \text{ and }
       \t\in\Gal\bigl(k(P)/k\bigr)
    \bigr\}.
\end{equation}
Clearly $S(P,n)\subset{E}\bigl(k(P)\bigr)$, since~$k(P)/k$ is Galois.

\claim{1}{$S(P,n)\subset E[n]\cap E\bigl(k(P)\bigr)$}
We have $nP\in{E}\bigl(k(nP)\bigr)$, and also~$k(nP)/k$ is Galois
since~$k(P)/k$ is abelian, from which it follows that
every~$\Gal\bigl(k(nP)/k\bigr)$-conjugate of~$nP$ is defined
over~$k(nP)$. Hence for all $\t\in\Gal\bigl(k(P)/k\bigr)$ and all
$\s\in\Gal\bigl(k(P)/k(nP)\bigr)$,
\[
  n(P^{\t\s}-P^\t) = (nP^\t)^\s - nP^\t = 0, 
\]
which shows that~$S(P,n)$ is contained in~$E[n]$.
The inclusion of~$S(P,n)$ in~$E\bigl(k(P)\bigr)$ follows
from the assumption that~$k(P)/k$ is Galois, so
every~$\Fbar/k$ conjugate of~$P$ is in~$E\bigl(k(P)\bigr)$.

\claim{2}{$S(P,n)$ is $\Gal\bigl(k(P)/k\bigr)$-invariant}
Let~$\t,\l\in\Gal\bigl(k(P)/k\bigr)$
and~$\s\in\Gal\bigl(k(P)/k(nP)\bigr)$.  Then
$\t\l\in\Gal\bigl(k(P)/k\bigr)$, and also
$\l^{-1}\s\l\in\Gal\bigl(k(P)/k(nP)\bigr)$, since~$k(nP)/k$ is Galois
(in fact, abelian).  Hence
\[
  (P^{\t\s}-P^\g)^\l = P^{\t\s\l}-P^{\t\l}
  = P^{\t\l(\l^{-1}\s\l)}-P^{\t\l} \in S(P,n).
\]
so~$\l$ maps~$S(P,n)$ to itself.

\claim{3}{$\bigl|S(P,n)\bigr| \ge \bigl[k(P):k(nP)\bigr]$}
The set~$S(P,n)$ contains in particular all points
of the form~\text{$P^\s-P$}
with~\text{$\s\in\Gal\bigl(k(P)/k(nP)\bigr)$} (i.e., take~$\t=1$),
and these points are distinct. Hence
\[
  \bigl|S(P,n)\bigr| \ge \bigl|\Gal\bigl(k(P)/k(nP)\bigr)\bigr| 
  = \bigl[k(P):k(nP)\bigr].
\]

We set the following notation, where note that Claim~2 tells us
that~$W$ is contained in~$V$:
\begin{align*}
  V&=E[n]\cong(\ZZ/n\ZZ)^2, \\
  \Aut(V)&=\Aut(E[n])\cong\GL_2(\ZZ/n\ZZ), \\
  W&=(\text{$\ZZ$-span of $S(P,n)$}) \subset V, \\
  \G(k) &= \Image\bigl(\Gal(\Fbar/k)\to\Aut(V)\bigr).
\end{align*}
Then we are in the following situation:
\begin{itemize}
\item
  $V$ is a free $\ZZ/n\ZZ$-module of rank~2.
\item
  $\G(k)$ is a subgroup of $\Aut(V)$.
\item
  $W$ is a $\G(k)$-invariant submodule of $V$ (from Claim 2).
\item
  The action of $\G(k)$ on $W$ is abelian (since
   $W\subset{E}(k(P))$ from Claim~1 and $k(P)/k$ is abelian).
\end{itemize}
These four conditions are exactly the assumptions needed to apply
Proposition~\ref{proposition:orderindexbd}, which yields the 
estimate
\begin{equation}
  \label{equation:WleIGki3}
  |W| \le \Index\bigl(\G(k)\bigr)^3.
\end{equation}
\par
The group~$\G(k)$ is the image of~$\Gal(\Fbar/k)$ in~$\Aut(V)$, but we
would like to replace it by the possibly larger group
\[
  \G(F) = \smash[b]{\Image\Bigl(\Gal(\Fbar/F)\to\Aut(V)\Bigr)}.
\]
Clearly~$\G(k)\subset\G(F)$, so
\begin{align}
  \label{equation:IGkiledIGF}
  \Index\bigl(\G(k)\bigr)
  &= \bigl(\Aut(V):\G(k)\bigr) \notag\\
  &=  \bigl(\Aut(V):\G(F)\bigr)\cdot\bigl(\G(F):\G(k)\bigr) \notag\\
  &\le  \bigl(\Aut(V):\G(F)\bigr)\cdot \bigl|\Gal(k/F) \bigr| \notag\\
  &\le  \bigl(\Aut(V):\G(F)\bigr)\cdot d 
     \qquad\text{since $[k:F]\le d$,} \notag\\
  &= d\cdot\Index\bigl(\G(F)\bigr).
\end{align}
Combining~\eqref{equation:WleIGki3} and~\eqref{equation:IGkiledIGF} yields
\begin{equation}
  \label{equation:Wled3IGF3}
  |W| \le d^3\cdot \Index\bigl(\G(F)\bigr)^3.
\end{equation}
\par
We now apply Serre's deep and fundamental theorem on the image of
Galois.

\begin{theorem}
\textup{(Serre~\cite{Serre72,Serre98})}
Let~$E/F$ be an elliptic curve defined over a number field.
For any prime~$\ell$, let
\[
  \rho_\ell : \Gal(\Fbar/F)\longrightarrow\Aut\bigl(T_\ell(E)\bigr)
\]
be the~$\ell$-adic representation attached to~$E/F$. 
Assume that~$E$ does not have complex multiplication. 
\begin{parts}
\Part{(a)}
The image of $\rho_\ell$ is of finite index in
$\Aut\bigl(T_\ell(E)\bigr)$ for all $\ell$.
\Part{(b)}
The image of $\rho_\ell$  equals 
$\Aut\bigl(T_\ell(E)\bigr)$ for all but finitely many~$\ell$.
\end{parts}
\end{theorem}

Serre's theorem is easily seen to be equivalent to the statement that
there exists a constant \text{$C_1(E/F)>0$} such that
\[
  \frac{\bigl|\Aut\bigl(E[N]\bigr)\bigr|}
  {\left|\Image\Bigl(\Gal(\Fbar/F)\to\Aut\bigl(E[N]\bigr)\Bigr)\right|}
  \le C_1(E/F)
  \qquad\text{for all $N\ge1$.}
\]
The crucial point here is that the constant~$C_1(E/F)$ is independent
of~$N$, so the index is bounded by a constant depending only on~$E/F$.
(Note that this is where we are using the assumption that~$E$ does not
have complex multiplication.) Applying Serre's estimate with~$N=n$
yields
\begin{equation}
  \label{equation:IGFleC1EF}
  \Index\bigl(\G(F)\bigr) \le C_1(E/F).
\end{equation}
\par
We combine the 
inequalities~\eqref{equation:Wled3IGF3}
and~\eqref{equation:IGFleC1EF} with Claim~3 
to obtain the estimate
\[
  \bigl[k(P):k(nP)\bigr]
  \le \bigl|S(P,n)\bigr|
  \le |W|
  \le  d^3\cdot \Index\bigl(\G(F)\bigr)^3
  \le  d^3\cdot C_1(E/F)^3.
\]
This this completes the proof of
Theorem~\ref{theorem:abstractindepresult}.
\end{proof}

%%%%%%%%%%%%%%%%%%%%%%%%%%%%%%%%%%%%%%%%%%%%%%%%%%%%%%%%%%%%%%%%%%%%%%%%
% Section: A direct sum decomposition via an idempotent relation
%%%%%%%%%%%%%%%%%%%%%%%%%%%%%%%%%%%%%%%%%%%%%%%%%%%%%%%%%%%%%%%%%%%%%%%%
\section{A direct sum decomposition via an idempotent relation}
\label{section:directsumdecomp}

Various versions of the results in this section are well-known, see
for example~\cite{Cornell91}
or~\cite[Theorem~6.3]{WalterKuroda79}.
%% *** Get the specific reference in Gary's article, he
%% indicates that the result on the prime-to-p part of the
%% class group for a (Z/pZ)^r extension is well known.
For the convenience of the reader we include proofs of the specific
statements that we require.

Let~$G$ be a  finite group. For each subgroup~$H\subset{G}$, the
associated \emph{idempotent}~$\e_H$ in the group ring of~$G$ is the
element
\[
  \e_H = \frac{1}{|H|}\sum_{\s\in H} \s \in \QQ[G].
\]
One easily verifies that~$\e_H^2=\e_H$. 

\begin{lemma}
\label{lemma:idempotent}
Let~$p$ be a prime, let~$G=\FF_p^r$, let~\text{$N=(p^{r-1}-1)/(p-1)$},
and let~$G_1,G_2,\ldots,G_N$ be the subgroups of~$G$ of index~$p$.
\par\vspace{2\jot}\noindent
\begin{tabular}{@{\enspace}l@{\quad}r@{}l}
\textup{(a)}
    &$\displaystyle\sum_{i=1}^N \e_{G_i}$
    &$\displaystyle{} = \frac{p^r-p}{p-1}\e_G + e$,\quad
      where $e\in G$ is the identity element.\\[4\jot]
\textup{(b)}
    &$\displaystyle\e_{G_i}\cdot\e_{G_j}$
    &$\displaystyle{} = \e_G$ \quad for all $i\ne j$. \\
\end{tabular}
\end{lemma}
\begin{proof}
Let
  $\hat{G}=\operatorname{Hom}(G,\FF_p)$
be the dual group to~$G$ and let~$\hat{G}^*$ denote the nonzero
elements of~$G$.  Then the kernel of each~$\chi\in\hat{G}^*$
is  an index~$p$ subgroup of~$G$, and~$\chi$ and~$\chi'$ give the same
subgroup if and only if~$\chi'=c\chi$ for some~$c\in\FF_p^*$. In other
words, the index~$p$ subgroups of~$G$ are in one-to-one correspondence with
the points of~$\hat{G}^*/\FF_p^*=\PP^{r-1}(\FF_p)$.
We also let~$G^*$ denote the nonzero elements of~$G$, i.e.,
\text{$G^*=\{\s\in{G}:\s\ne{e}\}$}, and compute
{\allowdisplaybreaks
\begin{align*}
  p^{r-1} \sum_{i=1}^N \e_{G_i}
  &= p^{r-1} \sum_{\chi\in\hat{G}^*/\FF_p^*} \e_{\Ker(\chi)} \\
%%  &=  \sum_{\chi\in\hat{G}^*/\FF_p^*} 
%%      \sum_{\substack{\s\in G\\\chi(\s)=0\\}} \s \\
  &=  \frac{1}{p-1}\sum_{\chi\in\hat{G}^*} 
      \sum_{\substack{\s\in G\\\chi(\s)=0\\}} \s \\
  &= \frac{1}{p-1}\sum_{\s\in G} 
           \bigl|\{\chi\in\hat{G}^* : \chi(\s)=0 \}\bigr| \s \\
  &= \frac{1}{p-1}  \left( |\hat{G}^*| e + 
         \sum_{\s\in G^*} 
                 \bigl|\{\chi\in\hat{G}^* : \chi(\s)=0 \}\bigr| \s \right) \\
  &= \frac{1}{p-1}\left( (p^r-1)e + (p^{r-1}-1)\sum_{\s\in G^*} \s \right) \\
  &= \frac{1}{p-1}\left((p^r-p^{r-1})e + (p^{r-1}-1)\sum_{\s\in G} \s \right)\\
  &= p^{r-1}e + \frac{p^{r-1}-1}{p-1}p^r\e_G.
\end{align*}
}%  This brace closes the \allowdisplaybreaks
This proves~(a).
\par
Next let~$i\ne{j}$ be distinct indices. Let~$\chi$ and~$\l$ be
elements of~$\hat{G}^*$ corresponding to~$G_i$ and~$G_j$, respectively.
Then
\[
  p^{2r-2}\e_{G_i}\cdot\e_{G_j}
   = \sum_{\s\in G_i}\sum_{\t\in G_j} \s\t 
   = \sum_{\substack{\s\in G\\ \chi(\s)=0\\}}
     \sum_{\substack{\t\in G\\ \l(\t)=0\\}} \s\t  
  = \sum_{\substack{\s\in G\\ \chi(\s)=0\\}}
     \sum_{\substack{\r\in G\\ \l(\s^{-1}\r)=0\\}} \r
  = \sum_{\rho\in G} \rho 
         \sum_{\substack{\s\in G\\ \chi(\s)=0\\ \l(\s)=\l(\rho)\\ }} 1.
\]
The fact that~$\chi$ and~$\l$ are distinct nonzero homomorphisms,
i.e., distinct nonzero linear maps, from~$\FF_p^r$ to~$\FF_p$
means that the pair of equations
\[
  \chi(\s) = c_1\qquad\text{and}\qquad \l(\s)=c_2
\]
has exactly~$p^{r-2}$ solutions for any given~$c_1,c_2\in\FF_p$.
(This is simply the number of points on the intersection of two transversal
hyperplanes in~$\AA^r(\FF_p)$.) Hence
\[
  p^{2r-2}\e_{G_i}\cdot\e_{G_j} = p^{r-2} \sum_{\rho\in G} \rho 
  = p^{2r-2} \e_G.
\]
Dividing by~$p^{2r-2}$ gives~(b).
\end{proof}

\begin{lemma}
\label{lemma:idempotentmodule}
Let~$p$,~$G$,~$N$, and~$G_1,\ldots,G_N$ be as in the
statement of Lemma~\ref{lemma:idempotent}. Let~$M$ be a
finite~$G$-module whose order is prime to~$p$, so in particular the
group ring~$\ZZ[p^{-1}][G]$ acts on~$M$.  
\begin{parts}
\Part{(a)}
$\e_HM=M^H$, i.e., $\e_HM$ is the subgroup of~$M$ fixed by~$H$.
\Part{(b)}
The  ``norm map''
\begin{equation}
  \label{equation:M=opluseGiM}
  \def\ds{\displaystyle}
  \begin{array}{ccc}
     \ds\frac{M}{M^G}
      &\ds\xrightarrow{\;\textup{norm}\;}
      &\ds\bigoplus_{i=1}^N \frac{\e_{G_i}M}{\e_GM}
       \cong\bigoplus_{i=1}^N \frac{M^{G_i}}{M^G} \\[4\jot]
      m  & \longmapsto & \ds \left(\e_{G_1}m,\ldots,\e_{G_N}m\right) \\
  \end{array}
\end{equation}
is an isomorphism. Its inverse is the summation map
\[
  (m_1,\ldots,m_N)  \longmapsto  m_1+\cdots+m_N.
\]
\end{parts}
\end{lemma}
\begin{proof}
First let~$\e_Hm\in\e_HM$ and let~$\t\in{H}$. Then
\[
  \t\e_Hm
  = \frac{1}{|H|} \sum_{\s\in H} \t\s m
  = \frac{1}{|H|} \sum_{\s\in H} \s m
  = \e_Hm,
\]
so~$\e_Hm\in M^H$. Conversely, if~$m\in M^H$, then
\[
  \e_H m
  = \frac{1}{|H|} \sum_{\s\in H} \s m
  = \frac{1}{|H|} \sum_{\s\in H}  m
  = m,
\]
so~$m=\e_H m\in\e_H M$. This proves that~$\e_HM=M^H$, which
completes the proof of~(a).
\par
Let~$\F$ be the map
\[
  \F:M
  \xrightarrow{\quad\textup{norm}\quad}
  \bigoplus_{i=1}^N \frac{\e_{G_i}M}{\e_GM},
  \qquad
  m \longmapsto \left(\e_{G_1}m,\ldots,\e_{G_N}m\right).
\]
It is clear that~$\e_GM$ is contained in the kernel of~$\F$.
Conversely, let~$m\in\Ker(\F)$. This means that there are~$m_i\in M$
satisfying \text{$\e_{G_i}m=\e_Gm_i$}. Summing over~$i$ and
using Lemma~\ref{lemma:idempotent}(a) yields
\[
  \e_G\sum_{i=1}^N m_i
  = \sum_{i=1}^N \e_{G_i} m
  = \left(\frac{p^r-1}{p-1}\e_G+e\right)m.
\]
Thus
\[
  m = \e_G\left(\frac{p^r-1}{p-1}m+\sum_{i=1}^r m_i\right)\in\e_GM,
\]
which gives the other inclusion. Hence~$\Ker(\F)=\e_GM$. 
\par
Next we show that~$\F$ is
surjective.  Let~$m_1,\ldots,m_N\in M$. We need to prove that the point
\[
  (\e_{G_1}m_1,\ldots,\e_{G_N}m_N) 
\]
is in the image of~$\F$. Let $m=\sum \e_{G_i}m_i$. Then for each~$j$
we use Lemma~~\ref{lemma:idempotent}(b)  to compute
\[
  \e_{G_j}m = \sum_{i=1}^N \e_{G_j}\e_{G_i}m_i
  = \e_{G_j}m_j  + \sum_{i\in j} \e_Gm_i
  \equiv \e_{G_j}m_j \pmod{\e_GM}.
\]
This proves that
\[
  \F\left(\e_{G_1}m_1+\cdots+\e_{G_N}m_N\right)
  = (\e_{G_1}m_1,\ldots,\e_{G_N}m_N) \pmod{\e_GM},
\]
completes the proof that~$\F$ induces an
isomorphism~\eqref{equation:M=opluseGiM} and that the inverse
isomorphism is the summation map.
\end{proof}

%%%%%%%%%%%%%%%%%%%%%%%%%%%%%%%%%%%%%%%%%%%%%%%%%%%%%%%%%%%%%%%%%%%%%%%%
% Section: An elementary Galois theory estimate
%%%%%%%%%%%%%%%%%%%%%%%%%%%%%%%%%%%%%%%%%%%%%%%%%%%%%%%%%%%%%%%%%%%%%%%%
\section{An elementary Galois theory estimate}
\label{section:galoistheoryestimate}

In this section we  prove a  basic estimate on the degree
of successive composita of Galois extensions,
cf{.}~\cite[Corollary~VI.1.15]{LangAlgebra}.

\begin{proposition}
\label{proposition:galdegest}
Let~$K_1,K_2,\dots,K_r$ be Galois extensions of a field~$k$, 
and for each \text{$2\le i\le r$}, let
\[
  K_i' = K_i \cap (K_1\cdots K_{i-1}).
\]
Then
\begin{equation}
  \label{equation:i1rKik}
  \prod_{i=1}^r [K_i:k] = [K_1\cdots K_r:k]  \prod_{i=2}^r [K'_i:k].
\end{equation}
\end{proposition}
\begin{proof}
Let~$E_1/k$ and~$E_2/k$ be any two Galois extensions
of~$k$. We claim that
\begin{equation}
  \label{equation:E1E2k}
  [E_1E_2:k]\cdot[E_1\cap E_2:k] = [E_1:k]\cdot[E_2:k].
\end{equation}
If $E_1\cap E_2=k$, then~\cite[Theorem~VI.1.14]{LangAlgebra} 
implies that~\eqref{equation:E1E2k} is true. The general
case of~\eqref{equation:E1E2k} can be reduced to this case as follows:
\begin{align*}
  [E_1E_2:k] &= [E_1E_2:E_1\cap E_2]\cdot[E_1\cap E_2:k] \\
  &= [E_1:E_1\cap E_2]\cdot[E_2:E_1\cap E_2]\cdot[E_1\cap E_2:k] \\
  &= \frac{[E_1:k]\cdot[E_2:k]}{[E_1\cap E_2:k]}.
\end{align*}
\par
We prove~\eqref{equation:i1rKik} by induction on~$r$. If~$r=1$, then
both sides are equal to~\text{$[K_1:k]$}. Assume now
that~\eqref{equation:i1rKik} is true for~$r$. We
use~\eqref{equation:E1E2k} with~$E_1=K_{r+1}$ and~$E_2=K_1\cdots{K}_r$
and note that \text{$E_1\cap E_2=K_{r+1}'$} to obtain
\[
  [K_{r+1}:k] =\frac{[K_1\cdots K_{r+1}:k]}{[K_1\cdots K_r:k]}
    \cdot [K_{r+1}':k].
\]
Multiplying~\eqref{equation:i1rKik} by this quantity
gives~\eqref{equation:i1rKik} with~\text{$r+1$} in place of~$r$, which
completes the induction.
%% To ease notation, let~$L_i=K_1\cdots{K}_i$, so 
%% \[
%%   L_{i+1} = L_iK_{i+1} \qquad\text{and}\qquad K_{i+1}' = K_{i+1}\cap L_i.
%% \]
%% Applying~\eqref{equation:E1E2k} with $E_1=L_i$ and $E_2=K_{i+1}$
%% yields
%% \[
%%   [L_{i+1}:k] [K_{i+1}' : k] = [L_i:k]\cdot [K_{i+1}:k].
%% \]
%% Now we multiply over $i=1,2,\ldots,r-1$, cancel common terms, and use
%% the fact that~\text{$L_1=K_1$} and~\text{$L_r=L$} to obtain the
%% desired result
\end{proof}

%%%%%%%%%%%%%%%%%%%%%%%%%%%%%%%%%%%%%%%%%%%%%%%%%%%%%%%%%%%%%%%%%%%%%%%%
% Section: Ideal class groups in multiquadratic fields
%%%%%%%%%%%%%%%%%%%%%%%%%%%%%%%%%%%%%%%%%%%%%%%%%%%%%%%%%%%%%%%%%%%%%%%%
\section{Ideal class groups in multiquadratic fields}
\label{section:idealclassgps}

For this section we fix the following notation.
\par\vspace{1\jot}\noindent
\renewcommand{\arraystretch}{1.1}
\begin{tabular}{cl}
$k/\QQ$&a Galois extension with group $\Gal(k/\QQ)=(\ZZ/2\ZZ)^r$.\\
$K$&the Hilbert class field of~$k$.\\
$N$&${}=2^{r-1}-1$.\\
$k_1,\ldots,k_N$&the distinct quadratic subfields of $k$.\\
$K_1,\ldots,K_N$&the Hilbert class fields of $k_1,\ldots,k_N$, respectively.\\
\end{tabular}
\renewcommand{\arraystretch}{1}
\par\vspace{1\jot}\noindent
Further,  for any finite abelian group~$G$, we
let~$G^\odd$ denote the largest subgroup of~$G$ of odd order,
and similarly~$h^\odd$ denotes the odd part of the integer~$h$.

\begin{proposition}
\label{proposition:idealclassgp2r}
With notation as set above, the natural restriction map
\[
  \Gal(K/k)^\odd\longrightarrow \prod_{i=1}^N \Gal(K_i/k_i)^\odd
\]
is an isomorphism.
\end{proposition}
\begin{proof}
Let~$H_k$ be the ideal class group of~$k$, and similarly for each~$i$,
let~$H_{k_i}$  denote the ideal class group of~$k_i$.  Standard
properties of the Artin map~\cite[X~\S1,~A2~\&~A4]{LangAlgNumbTh} give
a commutative diagram
\[\begin{CD}
  H_k @= H_k @>\text{Norm}>> H_{k_i} \\
  @V\text{Artin}V{\wr}V  @V\text{Artin}VV   @V\text{Artin}V{\wr}V  \\
  \Gal(K/k) @>>> \Gal(kK_i/k) @>\text{Res}>> \Gal(K_i/k_i) \\
  \end{CD}
\]
Combining these maps for~$i=1,2,\ldots,r$ and taking the odd parts of
each of the groups gives us a commutative diagram
\begin{equation}
  \label{equation:HkNHkiGKkGKiki}
  \begin{CD}
    H_k^\odd @>\operatorname{Norm}>> 
         \displaystyle\prod_{i=1}^N H_{k_i}^\odd \\
    @V\text{Artin}V{\wr}V   @V\text{Artin}V{\wr}V  \\
    \Gal(K/k)^\odd @>\operatorname{Res}>> 
	   \displaystyle\prod_{i=1}^N \Gal(K_i/k_i)^\odd\\
  \end{CD}
\end{equation}
\par
The Galois group~$\Gal(k/\QQ)=(\ZZ/2\ZZ)^r$ acts on~$H_k^\odd$, the
odd part of ideal class group of~$k$. Also note that the subgroups
of~$\Gal(k/\QQ)$ of index~$2$ are exactly the groups~$\Gal(k/k_i)$
for~\text{$1\le{i}\le{N}$}. This is exactly the situation needed
to apply Lemma~\ref{lemma:idempotentmodule}(b), which tells us
that there is an isomorphism
\begin{equation}
  \label{equation:HkHkGkQ}
  \frac{H_k^\odd}{(H_k^\odd)^{\Gal(k/\QQ)}}
  \xrightarrow[\sim]{\;\text{Norm}\;}
  \bigoplus_{i=1}^N 
  \frac{(H_k^\odd)^{\Gal(k/k_i)}}{(H_k^\odd)^{\Gal(k/\QQ)}}.
\end{equation}
(The norm maps $\Norm_{k/k_i}:H_k\to H_k^{\Gal(k/k_i)}$ appearing
in~\eqref{equation:HkHkGkQ} are actually the~$2^{N-1}$-power of the
``norm maps'' in Lemma~\ref{lemma:idempotentmodule}(b). However,
raising to the~$2^{N-1}$~power is an automorphism of~$H_k^\odd$, so it
is still valid to conclude that~\eqref{equation:HkHkGkQ} is an
isomorphism.)  \par In order to complete the proof of the proposition,
we use the following elementary lemma.

\begin{lemma}
\label{lemma:KFnHkHF}
Set the following quantities.
\par\vspace{1\jot}\noindent
\renewcommand{\arraystretch}{1.1}
\begin{tabular}{cl}
$k/F$&a Galois extension of number fields of degree~$n$.\\
$H_k'$&${}=H_k\otimes\ZZ[1/n]$,
        the prime-to-$n$ part of the class group of $k$.\\
$H_F'$&${}=H_F\otimes\ZZ[1/n]$,
        the prime-to-$n$ part of the class group of $F$.\\
\end{tabular}
\renewcommand{\arraystretch}{1}
\par\vspace{1\jot}\noindent
Then the natural map $H_F\to H_k$ induces an isomorphism
\[
  H_F'\to (H_k')^{\Gal(k/F)}.
\]
\end{lemma}
\begin{proof}
Consider the compositions
\begin{gather*}
  H_F'\longrightarrow (H_k')^{\Gal(k/F)}
   \xrightarrow{\Norm_{k/F}} H_F'\\
%%  \qquad\text{and}\qquad
\noalign{\noindent and}
  (H_k')^{\Gal(k/F)} \xrightarrow{\Norm_{k/F}} H_F' 
  \longrightarrow (H_k')^{\Gal(k/F)}.
\end{gather*}
Both compositions have the effect of raising to the
$n^{\text{th}}$~power.  Hence they are both isomorphisms, since~$H_F'$
and~$H_k'$ have order prime to~$n$.  Therefore each of
the individual arrows is also an isomorphism.
\end{proof}

We now resume the proof of
Proposition~\ref{proposition:idealclassgp2r}.  We apply
Lemma~\ref{lemma:KFnHkHF} with $F=k_i$ and $F=\QQ$ to obtain
\[
  (H_k^\odd)^{\Gal(k/k_i)} \cong H_{k_i}^\odd
  \qquad\text{and}\qquad
  (H_k^\odd)^{\Gal(k/\QQ)} \cong H_\QQ^\odd = 1.
\]
Substituting these values into the
isomorphism~\eqref{equation:HkHkGkQ} yields
\[
  H_k^\odd \xrightarrow[\sim]{\;\text{Norm}\;} \bigoplus_{i=1}^N H_{k_i}^\odd.
\]
Thus the top horizontal arrow in the commutative
diagram~\eqref{equation:HkNHkiGKkGKiki} is an isomorphism, and the
vertical arrows are Artin isomorphisms, hence the bottom arrow is also
an isomorphism.  This completes the proof of
Proposition~\ref{proposition:idealclassgp2r}.
\end{proof}

We now show that up to extensions of $2$-power
degree, the Hilbert class fields of distinct quadratic fields
are maximally disjoint.

\begin{proposition}
\label{proposition:KikiKjki2}
With the notation set at the beginning of this section,
for every~$2\le i\le r$ the degree
\[
  \biggl[ K_i\cap k_i\prod_{j<i}K_j \;:\; k_i \biggr]
\]
is a power of~$2$.
\end{proposition}
\begin{proof}
The diagram given in Figure~\ref{figure:fieldtower} should assist
in keeping track of the various fields under consideration.

%%%%%%%%%%%%%%%%%%%%%%%%%%%%%%%%%%%%%%%%%%%%%%%%%%%%%%%%%%%%%%%%%%%%%%
\begin{figure}
\begin{center}
\begin{picture}(260,295)(0,-10)
\thinlines
\put(0,0){\makebox(0,0){$\QQ$}}
\put(0,10){\line(0,1){20}}
\put(0,40){\makebox(0,0){$k_i$}}
\put(0,50){\line(0,1){20}}
\put(0,80){\makebox(0,0){$\displaystyle K_i\cap k_i\prod_{j<i}K_j$}}
\put(0,100){\line(0,1){40}}
\put(0,150){\makebox(0,0){$K_i$}}
\put(10,155){\line(3,1){110}}
\put(145,192){\makebox(0,0){$kK_i$}}
\put(145,185){\line(0,-1){40}}
\put(155,197){\line(3,1){70}}
\put(255,100){\line(0,1){115}}
\put(255,225){\makebox(0,0){$K_1\cdots K_N$}}
\put(255,235){\line(0,1){20}}
\put(255,265){\makebox(0,0){$K$}}
\put(10,40){\line(6,1){240}}
\put(255,90){\makebox(0,0){$k=k_1\cdots k_N$}}
\put(239,97){\line(-3,1){70}}
\put(40,85){\line(3,1){90}}
\put(145,125){\makebox(0,0){$\displaystyle kK_i\cap k\prod_{j<i}K_j$}}
\put(-5,20){\makebox(0,0){$\scriptstyle2$}}
\put(-5,60){\makebox(0,0){$\scriptstyle\textbf{?}$}}
\put(128,53){\makebox(0,0){$\scriptstyle2^{t-1}$}}
\put(80,90){\makebox(0,0){$\scriptstyle2^s$}}
\put(-40,95){\makebox(0,0)[c]
 {$\left\{\vrule height55pt depth55pt width0pt\right.$}}
\put(-55,95){\makebox(0,0){$H_i$}}
\put(292,177.5){\makebox(0,0)[c]
 {$\left.\vrule height87.5pt depth87.5pt width0pt\right\}$}}
\put(307,177.5){\makebox(0,0){$H$}}
\end{picture}
\end{center}
\caption{A tower of fields used in the proof of
Proposition~\ref{proposition:KikiKjki2}.  Various extensions of
$2$-power degree are indicated.  The proposition proves that the
extension marked with a~\textbf{?}  is also a~$2$-power extension.}
\label{figure:fieldtower}
\end{figure}
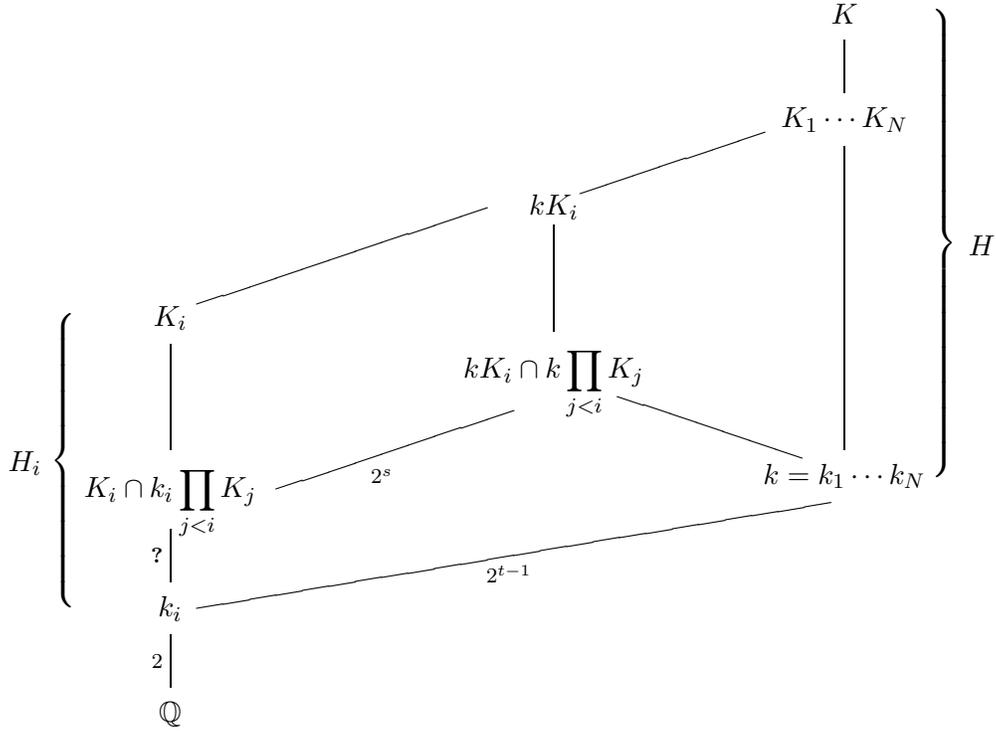
%%%%%%%%%%%%%%%%%%%%%%%%%%%%%%%%%%%%%%%%%%%%%%%%%%%%%%%%%%%%%%%%%%%%%%

We look first at the maps
\begin{equation}
  \label{equation:GKkGK1KNkGKiki}
  \Gal(K/k) \longonto \Gal(K_1\cdots K_N/k) 
  \longhookrightarrow \prod_{i=1}^N \Gal(K_i/k_i).
\end{equation}
Proposition~\ref{proposition:idealclassgp2r} tells us that if we
restrict to the odd parts of these groups, then the
composition~\eqref{equation:GKkGK1KNkGKiki} is an isomorphism. Since
the first map is surjective and the second map is injective, we conclude
that
\begin{equation}
  \label{equation:KkoddK1KNkoddKikiodd}
  [K:k]^\odd = [K_1\cdots K_N:k]^\odd = \prod_{i=1}^N [K_i:k_i]^\odd.
\end{equation}
\par
Next we apply Proposition~\ref{proposition:galdegest} to the base
field~$k$, the Galois extensions~$kK_1,\ldots,kK_N$ of~$k$, and 
their compositum \text{$L=K_1\cdots{K}_N$}. (Note that~$k$ is
already contained in~\text{$K_1\cdots{K}_N$}. For each
\text{$2\le{i}\le{N}$} we obtain the formula
\begin{equation}
  \label{equation:kKikkK1KNk}
  \prod_{i=1}^N [kK_i:k] = [K_1\cdots K_N:k]
  \prod_{i=2}^N
     \bigl[kK_i\cap(kK_1\cdots K_{i-1}) : k \bigr]
\end{equation}
Every extension~$k/k_i$ has degree a power of~$2$, so taking
the odd part of~\eqref{equation:kKikkK1KNk} and replacing~$k$
by~$k_i$ as appropriate yields
\begin{equation}
  \label{equation:KikioddK1KNkodd}
  \prod_{i=1}^N [K_i:k_i]^\odd = [K_1\cdots K_N:k]^\odd
  \prod_{i=2}^N
     \bigl[K_i\cap(k_iK_1\cdots K_{i-1}) : k_i \bigr]^\odd
\end{equation}
Comparing~\eqref{equation:KkoddK1KNkoddKikiodd}
and~\eqref{equation:KikioddK1KNkodd} gives the formula
\[
  \prod_{i=2}^N
     \bigl[K_i\cap(k_iK_1\cdots K_{i-1}) : k_i \bigr]^\odd = 1,
\]
and hence every factor in the product must equal~$1$.
\end{proof}

%%%%%%%%%%%%%%%%%%%%%%%%%%%%%%%%%%%%%%%%%%%%%%%%%%%%%%%%%%%%%%%%%%%%%%%%
% Section: The independence of Heegner points}
%%%%%%%%%%%%%%%%%%%%%%%%%%%%%%%%%%%%%%%%%%%%%%%%%%%%%%%%%%%%%%%%%%%%%%%%
\section{The independence of Heegner points}
\label{section:indepheegnerpts}

In this section we combine all of the previous material in order to
prove our principal result (Theorem~\ref{theorem:heegindepintro}),
which we restate here in a more precise form.

\begin{theorem}
\label{theorem:heegindep}
Let~$E/\QQ$ be an elliptic curve of conductor~$N$ and
let
\[
  \F_E:X_0(N)\longrightarrow E
\]
be a modular parametrization of~$E$.  There is a
constant~$C=C(E,\F_E)$ so that the following is true.
\par
Suppose that the following are given:
\par\vspace{1\jot}\noindent
\renewcommand{\arraystretch}{1.2}
\begin{tabular}{cl}
$k_1,\ldots,k_r$&\parbox[t]{.75\hsize}{\raggedright%
    distinct quadratic imaginary fields satisfying the
    Heegner condition for $N$.} \\
$h_1,\ldots,h_r$&the class numbers of $k_1,\ldots,k_r$.\\
$y_1,\ldots,y_r$&points  $y_i\in\CM^{(N)}(\Ocal_{k_i})$.\\
$P_1,\ldots,P_r$&the associated Heegner points, $P_i=\F_E(y_i)$.\\
\end{tabular}
\renewcommand{\arraystretch}{1}
\par\vspace{1\jot}\noindent
Assume further that the class numbers satisfy
\begin{equation}
  \label{equation:hioddgeC}
  h_i^\odd \ge C\qquad\text{for all $1\le i\le r$.}
\end{equation}
Then
\[
  \text{$P_1,\dots,P_r$ are independent in $E(\Qbar)/E_\tors$.}
\]
\end{theorem}
\begin{proof}
We will assume that~$P_1,\ldots,P_r$ are dependent and deduce an upper
bound for~$\min\{h_i^\odd\}$ that depends only on~$E$ and~$\F_E$.
Thus we assume that there is a relation
\begin{equation}
  \label{equation:n1P1nrPr=0}
  n_1P_1+\cdots+n_rP_r=0
  \qquad\text{with $n_i\in\ZZ$ not all~$0$.}
\end{equation}
Relabeling the points if necessary, we may assume that
\text{$n_r\ne0$}.  In order to complete the proof, it then suffices to
find a bound for~$h_r^\odd$
\par
We first apply Theorem~\ref{theorem:abstractindepresult} (with~$d=2$) to
deduce that
\begin{equation}
  \label{equation:krPrkrdivides}
  \bigl[k_r(P_r):k_r\bigr] 
  \enspace\text{divides}\enspace M\bigl[k_r(n_rP):k_r\bigr],
\end{equation}
where the constant~$M=M(E/F,2)$ is independent of~$P_r$ and~$n_r$.
\par
Next we let $k = k_1\cdots k_r$. Then~$\Gal(k/\QQ)=(\ZZ/2\ZZ)^t$ for
some~$t\le{r}$. We set~\text{$N=2^{t-1}-1$} and extend the list
of~$k_i$ to be the complete list~$k_1,\ldots,k_N$ of quadratic fields
contained in~$k$. Continuing with the notation from
Section~\ref{section:idealclassgps}, we set
\begin{align*}
  K&=\text{Hilbert class field of $k$,}\\
  K_i&=\text{Hilbert class field of $k_i$ for $1\le i\le N$.}
\end{align*}
\par
Each point~$P_i\in E(K_i)$, so the assumed linear
dependence~\eqref{equation:n1P1nrPr=0} tells us that
\[
  n_rP_r = -\sum_{i<r} n_iP_i \in E\Bigl(\prod_{i<r} K_i\Bigr).
\]
Thus
\[
  k_r(nP_r) \subset K_r \cap k_r \prod_{i<r} K_i,
\]
from which we conclude that~$\bigl[k_r(n_rP):k_r\bigr]$ divides
\begin{equation}
  \label{equation:KrkrKikr}
  \Bigl[ K_r \cap k_r \prod_{i<r} K_i : k_r \Bigr]
  =
  \Bigl[ K_r \cap k_r \prod_{i<r} K_i : k \Bigr]
  \cdot
  [k:k_r].
\end{equation}
Proposition~\ref{proposition:KikiKjki2} tells us that the first factor
on the righthand side of~\eqref{equation:KrkrKikr} is a power of~$2$, and
the second factor is clearly a factor of~$2$. This proves that
\[
  \text{$\bigl[k_r(n_rP):k_r\bigr]$ is a power of $2$.}
\]
Then~\eqref{equation:krPrkrdivides} tells us
that~$\bigl[k_r(P_r):k_r\bigr]^\odd$ divides~$M^\odd$, so it is
bounded by a constant depending only of~$E$.
\par
Finally, we use Proposition~\ref{proposition:heegnerptdegE}
and the fact that~$P_r$ is a Heegner point to deduce that
\[
  [K_r:k_r]^\odd \le (\deg\F_E)\bigl[k_r(P_r):k_r\bigr]^\odd
   \le (\deg\F_E)M^\odd.
\]
The quantity~$[K_r:k_r]^\odd$ is the odd part of the class number
of~$k_r$, so this contradicts~\eqref{equation:hioddgeC}.
Hence~$P_1,\ldots,P_r$ are independent.
\end{proof}

%%%%%%%%%%%%%%%%%%%%%%%%%%%%%%%%%%%%%%%%%%%%%%%%%%%%%%%%%%%%%%%%%%%%%%
%% *** Is the following remark worth making somewhere? 
%% \begin{remark}
%% There is a formula for the class number~$h_\Ocal$ of an order~$\Ocal$
%% in terms of the class number~$h$ of~$\Ocal_k$, 
%% the conductor~$c$ of~$\Ocal$,
%% and the behavior of the primes dividing~$c$ in~$k$. One version, given
%% in~\cite[Exercise~4.12]{Shimura71}, is
%% \[
%%   h_\Ocal = \#\Pic(\Ocal) = \frac{h\cdot c}{\#(\Ocal_k^*/\Ocal^*)}
%%      \prod_{p|c} \left[1-\genfrac{(}{)}{}{0}{k}{p}\frac{1}{p}\right].
%% \]
%% Here the Legendre symbol~$\genfrac{(}{)}{}{1}{k}{p}$ is defined to
%% be~$+1$,~$-1$, or~$0$, respectively, depending on whether~$p$ splits,
%% is inert, or ramifies in~$k$. Also note that
%% $\#(\Ocal_k^*/\Ocal^*)\le3$, and it is equal to~$1$ unless $k=\QQ(i)$
%% or~$\QQ(\rho)$.
%% \par
%% An equivalent alternative formula  given
%% in~\cite[Exercise~3.1]{DarmonCBMS04} is
%% \[
%%   h_\Ocal = \#\Pic(\Ocal) = 
%%    h\cdot \left|\frac{(\Ocal_k/c\Ocal_k)^*}{(\ZZ/c\ZZ)^*\Ocal_k^*}\right|.
%% \]
%% In any case, it's not hard to choose a value for~$c$ so that~$h_\Ocal^\odd$
%% is reasonably large.
%% \end{remark}
%%%%%%%%%%%%%%%%%%%%%%%%%%%%%%%%%%%%%%%%%%%%%%%%%%%%%%%%%%%%%%%%%%%%%%

%%%%%%%%%%%%%%%%%%%%%%%%%%%%%%%%%%%%%%%%%%%%%%%%%%%%%%%%%%%%%%%%%%%%%%%%
% Section: Deuring lifts, Heegner points, and the ECDLP
%%%%%%%%%%%%%%%%%%%%%%%%%%%%%%%%%%%%%%%%%%%%%%%%%%%%%%%%%%%%%%%%%%%%%%%%
\section{Deuring lifts, Heegner points, and the ECDLP}
\label{section:ecdlp}

In this section we briefly sketch an initially plausible approach to
solving the elliptic curve discrete logarithm problem using Heegner
points and explain why Theorem~\ref{theorem:heegindep} makes it
unlikely that this approach will yield anything better than an
algorithm with~$O(\sqrt{p})$ running time. We refer the reader
to~\cite{RosenSilverman05b} for further details.

\begin{definition}
An  \emph{Elliptic Curve Discrete Logarithm Problem} (ECDLP)
over~$\FF_p$ starts with a known elliptic curve~$E/\FF_p$ and
two points~$\Sbar,\Tbar\in{E}(\FF_p)$ and asks for the (smallest positive)
integer~$m$ such that
\[
    \Sbar = m \Tbar.
\]
(We use a bar to denote quantities defined over~$\FF_p$ or
to denote the reduction of a quantity modulo~$p$.)
It may be assumed that one knows
\[
  n_\Ebar = \#\Ebar(\FF_p).
%%  \qquad\text{and}\qquad
%%  a_\Ebar = p + 1 - n_\Ebar.
\]
\end{definition}

Suppose that it is possible to lift~$\Ebar$ to an elliptic curve~$E/\QQ$
of small conductor~$N$ and to find, in some reasonably explicit fashion,
a modular parametrization~\text{$\F_E:X_0(N)\to{E}$}. (If~$\Ebar$ is a
``random'' elliptic curve over~$\FF_p$, it is unlikely that this is
possible; but for reasons of efficiency in cryptographic applications,
it is not uncommon to take~$\Ebar$ to have very small coefficients.)
\par
A standard approach to solving the ECDLP is to choose many random
pairs of number~$(a_i,b_i)$ modulo~$n_\Ebar$ and try to find a
nontrivial relation amoung the points \text{$\Pbar_i=a_i\Sbar-b_i\Tbar$}.
If~$\sum{c_i}\Pbar_i=0$ is such a relation, then there is a good chance
that the resulting relation
\[
  \biggl(\sum_i c_ia_i\biggr)\Sbar = \biggl(\sum_i c_ib_i\biggr)\Tbar
\]
between~$\Sbar$ and~$\Tbar$ can be inverted to express~$\Sbar$ as a
multiple of~$\Tbar$. We thus look for a way of generating relations
among a given list~$\Pbar_1,\Pbar_2,\ldots$ of points
in~$\Ebar(\FF_p)$
\par
The map $\F_E:X_0(N)\to{E}$ has small degree, so there is a reasonable
chance that a randomly chosen point in~$E(\FF_p)$ will lift to
an~$\FF_p$-rational point on~$X_0(N)$. (The exact probability, which
we do not need, can be computed using the function field version of
the Tchebotarev density theorem.) Hence taking a subsequence, we
may assume that every point~$\Pbar_i\in{E}(\FF_p)$ lifts via~$\F_E$
to a point~$\ybar_i\in{X}_0(N)(\FF_p)$.
\par
In general, a point~$\ybar\in X_0(N)(\FF_p)$ corresponds to a
triple~$(\Abar,\Abar',\phibar)$ consisting of a pair of elliptic
curves~$\Abar/\FF_p$ and~$\Abar'/\FF_p$ and an isogeny
\[
  \phibar : \Abar\longrightarrow \Abar'
\]
defined over~$\FF_p$ with kernel $\Ker(\phibar)\cong\ZZ/N\ZZ$.  
Since~$\Abar$ and~$\Abar'$ are $\FF_p$-isogenous, they have the
same number of points, so we set the notation
\[
  n_\ybar=\#\Abar(\FF_p)=\#\Abar'(\FF_p)
  \qquad\text{and}\qquad
  a_\ybar=p+1-n_\ybar.
\]
(Note that these quantities can be computed
in polynomial time by the SEA variant of Schoof's
algorithm~\cite{Schoof85,Schoof95}.) Then the endomorphism rings of~$\Abar$
and~$\Abar'$ have discriminant
\[
  \D(\ybar) = a_\ybar^2 - 4p.
\]
Note that~$\D(\ybar)<0$ by Hasse's
theorem~\cite[Theorem~V.1.1]{SilvermanAEC86}.
\par
We perform this computation for each of the points~$y_1,y_2,\ldots$.
Taking a subsequence, we may assume that every integer~$\D(\ybar_i)$
is a fundamental discriminant. (A negative integer~$D$ is a
\emph{fundamental discriminant} if either it is odd, squarefree, and
$D\equiv1\pmod{4}$ or if it is divisible by~$4$ with~$D/4$ squarefree
and~$D/4\equiv2~\text{or}~3\pmod{4}$.)
This ensures that
\begin{multline*}
  \End(\Abar_i) = \End(\Abar_i') = \Ocal_{k_i}
  \quad\text{is the full ring of integers}\\
  \text{in the quadratic imaginary field }
  \smash[t]{k_i = \QQ\big(\sqrt{\D(\ybar_i)}\,\big)}.
\end{multline*}
In other words, both~$\Abar_i$ and~$\Abar_i'$ have~CM by the full ring
of integers of~$k_i$. (We also discard~$y_i$ in the unlikely event
that~$\Abar_i$ is supersingular, i.e., if~$a_{\ybar_i}=0$.)

\begin{remark}
In practice, we would like the fields~$k_1,k_2,\ldots$ to be
``not too independent,'' in the sense that we would like
the compositum~$k_1k_2\cdots$ to stop growing as we take more
and more points. This can be accomplished by keeping only those~$\ybar_i$
whose associated discriminant is $B$-smooth for an appropriately chosen
value of~$B$. For example, taking~$B=O(e^{\sqrt{\log p\log\log p}})$
as usual, suppose that the discriminants are reasonably randomly distributed
(a theorem of Birch~\cite{Birch68} says that they follow a
Sato-Tate distribution). Then
we can collect~$O(MB)$ points~$y_i$ and
fields~$k_i$ in time~$O(MB)$, while the compositum~$k_1k_2k_3\cdots$
is always contained in the field~$\QQ\bigl(\sqrt{\ell}:\ell\le{O}(B)\bigr)$,
independent of~$M$
\end{remark}

We next use a variant of Deuring's lifting theorem~\cite{Deuring41} to
lift~$\Abar$ and~$\Abar'$ to~CM elliptic curves defined over the
Hilbert class field~$K_i$ of~$k_i$ with the property that
\[
  \End(\Abar_i) = \End(A_i)\qquad\text{and}\qquad
  \End(\Abar_i')=\End(A_i).
\]
This can be done so that the cyclic
$N$-isogeny~\text{$\phibar:\Abar\to\Abar'$} lifts to a cylic
$N$-isogeny~\text{$\f:A\to{A'}$}.  (For modern expositions of
Deuring's theorem, see~\cite[13~\S5, Theorem~14]{LangEllipticFnc}
or~\cite{Oort73}).  The triple~$(A_i,A_i',\f)$ then corresponds to a
point~$y_i\in{X}_0(N)(K_i)$ that lifts~$\ybar_i$.
\par
Pushing these points forward, we obtain Heegner points
\[
  P_i = \F_E(y_i) \in E(K_i)
  \qquad\text{satisfying}\quad P_i\bmod{p} = \Pbar_i.
\]
(More precisely, there is a degree~1 prime ideal~$\gp_i$ in~$K_i$
so that~\text{$P_i\bmod{\gp_i}$} equals~$\Pbar_i$.)
\par
Our goal is to generate a list~$P_1,\ldots,P_r$ of points that are
dependent and to find an equation of dependency.  There are at least
two plausible methods to check whether $P_1,\ldots,P_r$ are dependent,
despite the fact that we generally cannot write down explicitly any
reasonable representation of their fields of
definition~$K_1,\ldots,K_r$.  (Note that
$[K_i:k_i]=h_{k_i}\approx\sqrt{|\Disc{k_i}|}$.)  First, we can try to
use the modular interpretation of the points~$y_i$ to compute the
canonical height pairing~$\langle{y_i},y_j\rangle$ in terms of other,
more easily computable, quantities. The algebraic and analytic
formulas developed by Gross, Kohnen and
Zagier~\cite{GrossKohnenZagier87} might be useful for this
approach. Second, we can use the standard proof of the weak
Mordell-Weil theorem to map a putative linear relation
from~$E(K_1\cdots{K}_r)$ into the Selmer group and attempt to derive
information about the coefficients of the relation. Note that there is
never any trouble checking if a potential relation is correct, since
it is easy to verify if~$\sum c_i\Pbar_i$ is zero, and if it is, then
we do not actually care if the relation lifts.
\par
Thus if the Deuring-Heegner lifts of the
points~$\Pbar_1,\Pbar_2,\ldots\in\Ebar(\FF_p)$ had a reasonable
probability of being dependent, then the algorithm outlined above
might yield a theoretical (and conceivably even a practical) algorithm
to solve the ECDLP on curves with small coefficients. However,
Theorem~\ref{theorem:heegindep} provides fairly convincing evidence
that this approach will not work, and indeed our initial motivation in
attempting to prove a theorem such as Theorem~\ref{theorem:heegindep}
was our desire to assess the effectiveness of Deuring-Heegner lifts
for solving the ECDLP

%%%%%%%%%%%%%%%%%%%%%%%%%%%%%%%%%%%%%%%%%%%%%%%%%%%%%%%%%%%%%%%%%%%%%%%%
% Section: Remarks on the $2$-part of the ideal class group
%%%%%%%%%%%%%%%%%%%%%%%%%%%%%%%%%%%%%%%%%%%%%%%%%%%%%%%%%%%%%%%%%%%%%%%%
\section{Remarks on the $2$-part of the ideal class group}
\label{section:2partclassgroup}

Let
\[
  \Dcal(X) = \{
    \text{fundamental discriminants $-D$ with $1\le D\le X$} \},
\]
and for each~$-D\in\Dcal(X)$, let~$h_D$ denote the class number of the
quadratic imaginary field~$\QQ(\sqrt{-D}\,)$.  
%% Further, for any positive integer~$h$, let~$h^\odd$ denote the odd
%% part of~$h$ as usual
In view of Theorem~\ref{theorem:heegindep}, it would be of
considerable interest to have some knowledge of the growth rate of the
counting function
\[
  N(C;X) = \#\bigl\{-D\in\Dcal(X) : h_D^\odd \le C \bigr\}.
\]
\par
Genus theory says that~$h_D$ is divisible by~$2^{\nu(D)}$,
where~$\nu(D)$ is the number of odd primes dividing~$D$. 
More precisely, genus theory tells us that the $2$-rank of the class
group is (essentially) equal to~$\nu(D)$. Since the~$n^{\text{th}}$~prime
is~$O(n\log{n})$, this immediately implies that
\[
  \text{($2$-rank of $H_D$)} \ll \frac{\log D}{\log\log D},
\]
so the $2$-part of~$h_D$ coming from the rank is negligible compared
to~$D$. However, little seems to be known about the expected growth of
the \text{$2$-exponent} of~$H_D$, see for
example~\cite{BoydKisilevsky72,Earnest89,WalterKuroda79}.
\par
If we look at all integers, then it is easy to see that 
\[
  \#\{n\le N : n^\odd\le C\}
  =  \left\lfloor \frac{C+1}{2}\right\rfloor\cdot\log_2(N) + O(1)
  \qquad\text{as $N\to\infty$.}
\]
The Brauer-Siegel theorem~\cite[Chapter~XVI,
Theorem~4]{LangAlgNumbTh} says that the magnitude of~$h_D$ is
approximately~$\sqrt{D}$, so~\text{$\{h_D:D\in\Dcal(X)\}$}
contains~$O(X)$ integers roughly less than~$\sqrt{X}$. This leads to the
following question, which we do not honor with the name of ``conjecture''
because there seems to be little evidence either for or against it.

\begin{question}
\label{question:countsmalloddpartclassnum}
Fix a constant~$C$. Is it true that
\[
  \#\bigl\{-D\in\Dcal(X) : h_D^\odd \le C \bigr\}
%%  N(C;X)
   \gg\ll \sqrt{X}\log(X)
  \qquad\text{as $X\to\infty$,}
\]
where the implied constants depend on~$C$?
\end{question}  

The numerical evidence is far from compelling. Indeed, the data in
Table~\ref{table:countoddclnumb} seems to suggest that~$N(C;X)\sim \k_C
X$ for some~$\k_C>0$, but it seems unlikely that this is true, and
indeed the value of~$\k_C$ shows a slow, but steady, decrease as we
eliminate the smaller data from the top of the table. However, the
data also does not suggest that~$N(C;X)\ll{X}^{1-\d}$ for any
particular~$\d>0$.

\begin{table}
\begin{center}
%%\small
\caption{Counting quadratic imaginary fields whose class number has
small odd part and whose discriminant~$-D$ satisfies
\text{$500{,}000\le{D}\le1{,}000{,}000$} (computations performed with PARI-GP).}
\label{table:countoddclnumb}
\tiny
\begin{tabular}{|r|r|r|r|r|r|} \hline
\# of $D$ &  $h^\odd\le1$ &  $h^\odd\le3$ &  $h^\odd\le5$ &
 $h^\odd\le7$ &  $h^\odd\le9$ \\
%% \# of $D$ & 1 & 3 & 5 & 7 & 9 \\
\hline \hline
5555 & 134 & 312 & 479 & 617 & 770 \\
11110 & 264 & 589 & 901 & 1183 & 1486 \\
16665 & 407 & 899 & 1381 & 1793 & 2273 \\
22220 & 525 & 1202 & 1825 & 2375 & 3012 \\
27775 & 667 & 1511 & 2281 & 2954 & 3730 \\
\hline
33330 & 775 & 1789 & 2699 & 3497 & 4427 \\
38885 & 893 & 2062 & 3114 & 4052 & 5154 \\
44440 & 1037 & 2358 & 3544 & 4599 & 5847 \\
49995 & 1163 & 2641 & 3970 & 5140 & 6532 \\
55550 & 1278 & 2909 & 4374 & 5681 & 7225 \\
\hline
61105 & 1395 & 3190 & 4786 & 6219 & 7902 \\
66660 & 1515 & 3479 & 5216 & 6785 & 8610 \\
72215 & 1622 & 3745 & 5612 & 7302 & 9260 \\
77770 & 1735 & 3995 & 6003 & 7816 & 9922 \\
83325 & 1862 & 4272 & 6406 & 8343 & 10594 \\
\hline
88880 & 1983 & 4557 & 6820 & 8887 & 11278 \\
94435 & 2093 & 4799 & 7173 & 9361 & 11891 \\
99990 & 2208 & 5073 & 7566 & 9867 & 12560 \\
105545 & 2310 & 5295 & 7938 & 10344 & 13170 \\
111100 & 2427 & 5556 & 8348 & 10875 & 13813 \\
\hline
116655 & 2533 & 5789 & 8714 & 11347 & 14419 \\
122210 & 2652 & 6039 & 9085 & 11831 & 15060 \\
127765 & 2748 & 6297 & 9465 & 12337 & 15697 \\
133320 & 2847 & 6541 & 9819 & 12799 & 16298 \\
138875 & 2946 & 6789 & 10184 & 13284 & 16938 \\
\hline
144430 & 3049 & 7012 & 10545 & 13753 & 17529 \\
149985 & 3142 & 7235 & 10880 & 14211 & 18127 \\
\hline\hline
\begin{tabular}{@{}c@{}}Linear\\correlation\\\end{tabular}
  &99.926\% & 99.932\% & 99.943\% & 99.949\% & 99.955\% \\
\hline
\end{tabular}
\end{center}
\end{table}

%%%%%%%%%%%%%%%%%%%%%%%%%%%%%%%%%%%%%%%%%%%%%%%%%%%%%%%%%%%%%%%%%%%%%%%%
% Bibliography
%%%%%%%%%%%%%%%%%%%%%%%%%%%%%%%%%%%%%%%%%%%%%%%%%%%%%%%%%%%%%%%%%%%%%%%%

\providecommand{\bysame}{\leavevmode\hbox to3em{\hrulefill}\thinspace}
\providecommand{\MR}{\relax\ifhmode\unskip\space\fi MR}

%% \bibliographystyle{amsalpha}
%% \bibliography{ECDLP}

\end{document}